\documentclass[12pt,reqno]{article}
\usepackage[T2A]{fontenc}
\usepackage{amsmath,amsthm,amssymb,amsfonts}
\usepackage{verbatim}

\numberwithin{equation}{section}
\usepackage[linktocpage=true]{hyperref}

\usepackage{tikz-cd}
\usepackage[all]{xy}
\usepackage{mathtext}
\usepackage[T1,T2A]{fontenc}
\usepackage[utf8]{inputenc}
\usepackage[russian,english]{babel}
\usepackage{amsfonts,longtable,amssymb,amsmath,latexsym, dsfont}
\usepackage{wrapfig}

\usepackage{mathtools}
\usepackage{wrapfig}
\usepackage{csquotes}
\usepackage{titlesec}
\titleformat{\subsubsection}[runin]
 {\normalfont\bfseries}
 {\thesubsubsection}
 {0.5em}
 {}
 [.]

\usepackage{amsthm, amsmath, amssymb, amsfonts, amscd, amscd, geometry}
\newsavebox{\ssa}
\def\a{\alpha}

\def\s{\sigma}

\def\G{\Gamma}

\def\n{\nabla}
\renewcommand{\gg}{\mathfrak{g}}

\newcommand{\R}{\mathbb R}
\newcommand{\N}{\mathbb N}

\newcommand{\Cl}{\text{Cl}}
\newcommand{\Z}{\mathbb Z}
\renewcommand{\C}{\mathbb{C}}

\theoremstyle{definition}
\newtheorem{theorem}{Theorem}[section]
\newtheorem{Th}[theorem]{Theorem}
\newtheorem{Prop}[theorem]{Proposition}
\newtheorem{Cor}[theorem]{Corollary}
\newtheorem{Lem}[theorem]{Lemma}
\newtheorem{Def}[theorem]{Definition}

\newtheorem{Remark}[theorem]{Remark}
\newtheorem{Example}[theorem]{Example}
\def\Cl{\text{\upshape Cl}}
\def\bt{\begin{Th}}
	\def\et{\end{Th}}
\def\bp{\begin{Prop}}
	\def\ep{\end{Prop}}
\def\bc{\begin{Cor}}
	\def\ec{\end{Cor}}
\def\bl{\begin{Lem}}
	\def\el{\end{Lem}}
\def\bd{\begin{Def}}
	\def\ed{\end{Def}}
\def\be{\begin{equation}}
	\def\ee{\end{equation}}
\def\bR{{\mathbb R}}
\def\bZ{{\mathbb Z}}
\def\bC{{\mathbb C}}
\def\bA{{\mathbb A}}
\def\bH{{\mathbb H}}

\def\bO{{\mathbb O}}
\def\bK{{\mathbb R}}
\def\bT{{\mathbb T}}

\def\cg{\mbox{\mcal G}}

\def\cn{\mbox{\mcal N}}

\def\R{\mbox{\mcal R}}

\def\ct{\mbox{\mcal T}}

\def\cT{\mbox{\mcal T}}
\def\cH{\mbox{\mcal H}}
\def\cN{\mbox{\mcal N}}
\newfont{\mcal}{eusm10 scaled \magstep1}

\def\cg{\mbox{\mcal G}}

\def\cn{\mbox{\mcal N}}

\def\Id{\mathrm{Id\;}}
\def\tr{\mathrm{tr\;}}

\def\Hom{\mathrm{Hom\;}}

\def\ad{\mathrm{ad}}

\def\GL#1{{\mathrm{GL(#1)}}}
\usepackage{tikz}
\usetikzlibrary{arrows}
\usepackage{tkz-berge}
\usepackage{tkz-graph}

\def\arr{\begin{array}{rlll}}
	\def\ea{\end{array}}
\def\be{\begin{equation}}
\def\ee{\end{equation}}

\def\bea{\begin{eqnarray}}
	\def\eea{\end{eqnarray}}
\def\bean{\begin{eqnarray*}}
	\def\eean{\end{eqnarray*}}
\newcommand{\Nil}{\operatorname{Nil}}

\newcommand{\Herm}{\operatorname{ Herm}}

\theoremstyle{remark}
\newtheorem{rem}{Remark}[section]
\usepackage{indentfirst}
\usepackage[style=alphabetic,citestyle=alphabetic,maxnames=99,maxcitenames=99, maxalphanames=99, sorting=nyt]{biblatex}

\begin{document}
	\begin{center}
		{\LARGE\bf
			Special Vinberg cones of rank 4.}
		\medskip
		\medskip
		
		D.V. Alekseevsky\footnote{Higher School of Modern Mathematics MIPT, Russian Federation.}  P. Osipov\footnote{National Research University Higher School of Economics, Russian Federation.}\footnote{Pavel Osipov was supported by the Basic Research Program of the National Research University Higher
School of Economics.}
	\end{center}

	\tableofcontents
	
	\newpage 
\begin{abstract}
 E.B. Vinberg developed a theory of homogeneous convex cones $C \subset V= \bR^n$, which has many applications. He gave a construction of such cones in terms of non-associative rank $n$ matrix T-algebras $\ct$, that consist of vector-valued $n \times n$ matrices $X = ||x_{ij}||, \, x_{ij} \in V_{ij} $ where $V_{ij}$ are Euclidean vector spaces. The multiplication in a T-algebra is determined by a system of isometric maps $V_{ij} \times V_{jk} \to V_{ik}$, s.t. $|v_{ij}\cdot v_{jk}| = |v_{ij}|\cdot |v_{jk}|$ that satisfies some axioms. A T-algebra is determined by its {\bf associative } subalgebra of upper triangular matrices $\cg$ or its niladical $\cn$, called the Nil-algebra. The connected Lie group $G\subset \cg$ of the upper triangular (non-degenerate) matrices acts in the vector space $Herm_n \subset\ct$ of Hermitian matrices and the orbit $C = G(I)\subset Herm_n$ of the identity matrix $I$ is a convex cone with a simply transitive action of $G$. Conversely, any homogeneous convex cone is obtained by this construction. 
 
Generalizing the notion of rank 3 Clifford T-algebra \cite{alekseevsky2021special}, \cite{alekseevsky2021special1}, we define notions of rank $n$ special T-algebra and Clifford Nil-algebra, which define a special Vinberg cone. 
We associate with a Clifford Nil-algebra $\cN$ a directed acyclic graph $\G=\G(\cN)$ of diameter 1 and show that Clifford Nil-algebras with given graph $\G$ bijectively correspond to its admissible equipments. This gives an effective method of classification of Clifford Nil-algebras and associated special Vinberg cones. We apply this approach for explicit classification of rank 4 special Vinberg cones.
 
\end{abstract}

	\section{Introduction}	

 \subsection{Vinberg's theory of homogeneous convex cones}
 \subsubsection{Basic definitions}
A {\bf convex cone} is an open convex $\bR^{>0}$-invariant domain $C$ without straight lines in a vector space $V=\bR^n$. It is called {\bf homogeneous}, if the linear group
 $$ G= \{ A \in \mathrm{GL}(V), A(C) = C \}
 $$
 of automorphisms acts transitively in $C$. Then there exists a solvable subgroup $G \subset \mathrm{Aut}(C) $ (the Vinber group)
 which acts {\bf simply transitively} in $C$ .
 The {\bf dual cone }
 $C^* \subset V^*$ in the dual space is defined by 
 $$C^* = \{\xi \in V^*,\, \xi(V)>0 \, \forall X \in C \}.$$
 A homogeneous convex cone is called {\bf self-dual} if $g\circ C = C^*$ for some Euclidean metric $g$. 
 Then the dual cone $C^*$ is also convex and homogeneous.
 
 A homogeneous convex cone $C \subset V$ admits a unique up to scaling $\mathrm{Aut}(C)$
invariant measure $ \chi(x) dx := \chi(x) dx^1 \wedge \cdots \wedge dx^n$, where 
$$\chi (x) = \int_{C^*}e^{-\langle\xi,x\rangle}d\xi$$
is the density i.e.
$$\chi(AX) = \det A^{-1}\chi(x),\, A \in \text{Aut}(C).$$
 It is called the {\bf characteristic function of the cone} and it 
completely characterizes the cone $C$. 

 The 1-form $-d \chi$ called the {\bf Koszul form} and the map $ \star: =-d \chi: C \to C^* $ defines an $\text{Aut}(C)$-equivariant diffeomorphism of the cone $C$ onto the dual cone $C^* \subset V^*$. Moreover, if $C$ is homogeneous, then $\star^2 = \Id $.
The {\bf canonical } $Aut(C)$-invariant Riemannian Hessian metric is defined by
$$g _C = \text{Hess} \log \chi = \partial^2 \log \chi_C.$$ 

The theory of homogeneous convex cones was developed by E.B.Vinberg in the seminal paper \cite{vinberg1963theory}. This theory and,
in particular, the characteristic function have many applications in differential geometry, physics, convex programming, information geometry, statistics, etc.
The above Vinberg paper is the most cited work in his rich heritage.
 
 Originally, the problem of classification of homogeneous convex cones emerged from Cartan's problem on classification of homogeneous under holomorphic automorphisms bounded domains in $\bC^n$ (\cite{пятецкий1957оценке}, \cite{пятецкий1961геометрия}). If $C\subset \bR^n$ is a homogeneous convex cone, then the domain $D= \bR^n \oplus iC$ is biholomorphic to a homogeneous bounded domain in $\bC^n$. Homogeneous bounded domains, obtained by this construction, are called Siegel domains of the first kind. There exist also Siegel domains of the second kind. Any homogeneous bounded domain in $\bC^n$ is a Siegel domain of the first or second kind. 
 
\subsubsection{Vinberg theory of T-algebras and homogeneous convex cones}
 
In the paper \cite{vinberg1963theory} E. B. Vinberg proved that any homogeneous convex cone can be described as the positive cone $C$ in the space $\mathrm{Herm}_n \subset \mathcal{T}$ of Hermitian matrices of some matrix T-algebra $\mathcal{T}$. 
The positive cone $C = G(\Id)$ is defined as the $G$-orbit of the identity matrix, where $G \subset \mathcal{T}$ is the Vinberg group of upper triangular matrices with positive diagonal elements.\\
The matrix T-algebra of rank $n$
$$(\mathcal{T} =\sum_{i,j=1 }^n V_{ij}, g(\cdot, \cdot) = <\cdot,\cdot>)$$ is an orthogonally bigraded Euclidean vector space, such that $V_{ii} = \bR, V_{ji} = V_{ij}^*, \forall i \neq j$, with bilinear multiplication defined by a system of bilinear isometric maps 
$$
\mu_{ijk} : V_{ij} \times V_{jk} \to V_{ik}, \ \ \ (v_{ij }, v_{jk}) \to v_{ij} \cdot v_{jk}
$$
such that 
$$ \langle v_{ij } \cdot v_{jk}\rangle^2 = \langle v_{ij}, v_{ij}\rangle \cdot \langle v_{jk}, v_{jk}\rangle. $$
It is convenient to consider elements $X= || x_{ij}|| \in \mathcal{T}$ as matrices with a non-associative matrix multiplication, defined by the isometric maps $\mu_{ijk}$.

The matrix multiplication must satisfy some system of axioms. The most
important one is that the subalgebra $\cN= \cN (\mathcal{T})$ of the strictly upper triangular matrices is an associative algebra. In fact, this condition (with a weak additional condition) is sufficient to reconstruct the T-algebra $ \mathcal{T}$ and the associated homogeneous convex cone.

Vinberg’s main theorem states that any homogeneous convex cone is realized as the positive cone 
$$C = G(\Id)= \{ AA^*,\, A \in G\} \subset Herm_n \subset \mathcal{T}$$
of Hermitian matrices of some T-algebra, that is, the orbit of the identity matrix $\Id$ under the action of the Vinberg upper triangular group $G$, generated by the Lie algebra $\gg, $ with the Lie bracket $[A,B] = AB - BA$.

 The dual T algebra $\mathcal{T}'$ is obtained from $\mathcal{T}$
by the symmetry $\tau$ across the antidiagonal.
The algebras $\mathcal{T},\mathcal{T}' $ are isomorphic as abstract algebras. They are isomorphic as matrix algebras if and only if the cone $C$ is self-dual.
According to K\"ocher-Vinberg theorem, this is the case if the cone $C$ with the canonical metric is a Riemannian symmetric space, more precisely, the direct product of the cones of Hermitian positively definite matrices over division algebra $\bR,\bC,\bH, $ and for rank $3$ over $\mathbb{O}$.
 
\subsection{Special Vinberg cones of rank 3}
 The first problem for an explicit description of homogeneous convex cones is to specify the system of isometric maps $\mu: X \times Y \to Z$ between Euclidean spaces. If $\dim Y =\dim Z$ then the isometric map $\mu$ defines the structure of a $\bZ/2\bZ$-graded Clifford module in the space $S = S^0\oplus S^1 = Y \oplus Z $ over $V=X$
 and the Clifford multiplication 
 $$ 
 \mu : V \times S^0 \to S^1, \ \ \ \mu' : V \times S^1 \to S^0
 $$
 defines the isometric map 
 $\mu : X \times Y \to Z,\, (x, y) \to \mu_x y = x \cdot y $ 
 and the conjugated isometric maps $\mu' : X \times Z\to Y $.\\
 We call such isometric maps $\mu, \mu'$ {\bf Clifford} (isometric) maps.\\
In \cite{alekseevsky2021special}, the authors defined and studied a class of rank 3 T-algebras $\mathcal{T}$ (called {\bf special}) and the corresponding homogeneous convex cones $C$, (called {\bf special Vinberg cones}) associated to a $\bZ/2\bZ$-graded Clifford module $S = S^0 +S^1$.
 More precisely, any Clifford $Cl(V)$ 
graded module $S = S^0+S^1$ defines two dual T-algebras
$\mathcal{T} = \mathcal{T}(S) $ and , resp., $\mathcal{T}'=\mathcal{T}'(S)$, consisting of matrices of the form
$$X =||X_{ij}||=
\begin{pmatrix} 
x_1&v&s_1\\
v'&x_2& s_0 \\
 s'_1 & s'_0& x_3
\end{pmatrix} \ \
\text {and, resp.}, \ \
X^\tau =||X^\tau_{ij}||=
\begin{pmatrix} x_1&s_0&s_1\\
s'_0&x_2& v \\
s'_1 & v' & x_3 
\end{pmatrix}.
$$
where $\tau$ is the anti-transposition and $x_i \in \bR, v \in V,v' \in V^*, s_0 \in S^0, s'_0 \in (S^0)^*, s_1 \in (S^1)^*$.
The Vinberg group $G $ of the cone $\mathcal{T}$ consists of the matrices of the form 
$$ A = ||a_{ij}||=
\begin{pmatrix}
 a_{11}&a_{12}& a_{13}\\
 0&a_{22}& a_{23}\\
 0&0& a_{33}
\end{pmatrix}, a_{ii}\in \bR^+,\, a_{12} \in V, \, a_{13} \in S_1, \, a_{23} \in S_0.
$$ 
and similarly for the Vinberg group $G'$ of the dual T-algebra.
 This construction also works in the pseudo-Euclidean case and associates with a $\bZ/2\bZ$-graded Clifford $Cl(V)$-module $S=S^0+S^1$ over a pseudo-Euclidean space $V$ a homogeneous open, but not necessarily convex, cone.
 
 Since a Vinberg upper triangular group $G$ acts simply transitively in the Vinberg cone, the coordinates of the matrix elements $a_{ij}, \, i\leq j$ of a matrix $A \in G$ define the coordinates of the cone, called group coordinates. 
 In \cite{alekseevsky2021special}, the authors explicitly calculated the characteristic function $\chi$ of a Vinberg rank 3 cone $C(S)$ and three homogeneous polynomials $p_3 (X) = X_{33} , p_2(X), p_1(X)$ in $\mathcal{T}(S)$. According to Vinberg, these polynomials determine the cone $C(S)$ by inequalities $p_1(X) >0, \, p_2(X) >0, \, p_3(X) >0$. They have degree 1,2 and 4 and their restriction to the cone in a point $X = A A^*$ in group coordinates $a_{ij}$ are given by
$$p_3(AA^*)= a_{33}^2,\, p_{2}(AA^*)= (a_{22}a_{33})^2,\, 
p_1(AA^*)= a_{11}^2a_{22}^2a_{33}^4. $$

 The function $d(X) := \frac{p_1(X)}{p_3(X)}$ on a Vinberg rank 3 T-algebra is called the {\bf determinant function} 
since its restriction $d(X)|_C = d(AA^*)$ to the cone in group coordinates is given by 
 $$ d(AA^*) = det(A )^2: = (\a_{11} a_{22} a_{33})^2.$$
 The function $d(X)$ for the dual cone $C^*(S)$ is a homogeneous cubic and it defines a 
 special real manifold in the sense of de Wit and Van Proeyen \cite{de1992special}).
 Surprisingly, the determinant function of the Vinberg cone $C(S)$ is only a homogeneous rational function of degree three. The papers \cite{alekseevsky2021special}, \cite{alekseevsky2021special1} contain the erroneous statement that for the cone $C(S)$ the function $d(X)$ is also a cubic polynomial. It is corrected in \cite{alekseevsky2023erratum} and \cite{matrix}. 

 \subsection{The main results}
We define a special class of T-algebra, called {\bf Clifford T-algebras}, as T-algebras with multiplication, defined by Clifford isometric map. The associated homogeneous convex cone is called a {\bf (special) Vinberg cone}. 
 
\subsubsection{Nil-algebras} A Vinberg result implies that a special T-algebra $\mathcal{T}$ is completely reconstructed from its associative subalgebra $\gg$ of the upper triangular matrices and even from its nilpotent subalgebra, called {\bf Nil-algebra}. Two Nil-algebras are called equivalent if the corresponding cones are isomorphic. \\
A Nil-algebra is intrinsically characterized as an associative nilpotent algebra of upper triangular matrices
$
\cN = \{ A = \{ a_{ij}\} \ | \ \forall i<j \ : \ a_{ji}=0, \ a_{ii}=0, \ a_{ij} \in \cN_{ij} \} 
$
with the associative multiplications $\mu_{ijk} : \cN_{ij }\times \cN_{jk}\to \cN_{ik}$ defined by isometric maps and satisfying the Vinberg condition \eqref{Vinberg condition}. \\
It is called {\bf Clifford Nil-algebra} if for all triples $i<j<k$ the isometric map $\mu_{ijk}$ is Clifford, that is $\mu_{ijk}$ 
 coincides with the multiplication map $V\times S^0 \to S^1$ or $S^0\times V \to S^1$, where $V$ is a Euclidean space and $S^0\oplus S^1$ is a $\bZ/2\bZ$-graded $\Cl(V)$-module. Two Nil-algebras $\cN$ and $\cN'$ of rank $n$ are equivalent if there exists an isomorphism of algebras $f: \cN \to \cN'$ and a permutation $\sigma \in S_n$ such that for any $1\le i\le j \le n$ we have $f(\cN_{ij})=\cN_{\sigma(i),\sigma(j)}$.

 The T-algebra $\mcal{T}$, defined by a Clifford Nil-algebra $\cN$ is called a {\bf Clifford T-algebra }
 and the corresponding homogeneous convex cone $C \subset \text{Herm}_n \subset \mcal(T)$ is called {\bf (special) Vinberg cone} of rank $n$.
 
 The classification of special Vinberg cones reduces to the classification of Clifford Nil-algebras.
 
\subsubsection{The main results}
We propose a method for the description of special Vinberg cones of rank $n$ in terms of associated Nil-algebras and apply it to the classification of Vinberg cones of rank 4. 

Removing the associativity condition and the Vinberg condition, we define the notion of {\bf quasiNil-algebra} and {\bf Clifford quasiNil-algebras}. Then we associate with any indecomposable quasiNil-algebra (that is not equivalent to a direct sum of two other such algebras) $\cN$ a directed acyclic graph $\G = \G(\cN)$ of diameter $d(\G)=1$ called the {\bf adjacency graph} of the quasiNil-algebra $\cN$. 
A directed acyclic graph of diameter 1 is called a Nil-graph. Any Nil-graph $\G$ is the 
 adjacency graph $\G = \G(\cN)$ of some indecomposable Nil-algebra $\cN$ and all quasiNil-algebras with given Nil-graph $\G$ bijectively correspond to the equipments of $\G$.
 
The description of equipments of the graph $\G$, hence Clifford quasiNil-algebras is a pure combinatorial problem.

 We solve it for all Nil-graphs with $n=4$ vertices and get a classification of rank 4 Clifford quasiNil-algebras. Then the classification of rank 4 Nil-algebras reduces to the description of when such quasiNil-algebras are associative and satisfy the Vinberg condition \eqref{Vinberg condition}. 

 The main result is the classification of rank 4 special Vinberg cones in terms of associated rank 4 Nil-algebras. We say that a quasiNil-algebra $\cN$ has the {\bfseries nilpotency index} $NI(\cN)=m$ if the maximal number of elements with nonzero products is $m-1$, or , equivalently, the maximal length $ML(\G) $ of a path in the adjacency graph $\G = \G(\cN) $ is $ML(\G) =m$.
 
 In the main theorem below, we list indecomposable matrix Nil-algebras $\cN$ of rank $n=4$ up to duality (defined by the antitransposition). The adjacency graphs of the dual-to-$\cN$ quasiNil-algebra have the dual-to-$\G(\mcal{N})$ graph, which is obtained from $\G(\mcal{N})$ by reversing arrows.
 
 \begin{theorem}\label{index2} \ 
 Below is the list of rank 4 Clifford Nil-algebras of $ NI \le 2$ up to a duality.
 
 \begin{itemize} 
 \item [a)] Case $NI =1$
		$$
		\begin{pmatrix}
			0 &V_{1}& V_2 & V_3 \\
			0 & 0 & 0 & 0 \\
			0 & 0 & 0 & 0 \\
			0 & 0 & 0 & 0
		\end{pmatrix}, \ \ \ 
		\begin{pmatrix}
			0 & 0 & V_1 & V_2 \\
			0 & 0 & 0 & V_3 \\
			0 & 0 & 0 & 0 \\
			0 & 0 & 0 & 0
		\end{pmatrix}, 
		\ \ \ 
		\begin{pmatrix}
			0 & 0 & V_1 & V_2 \\
			0 & 0 & V_3 & V_4 \\
			0 & 0 & 0 & 0 \\
			0 & 0 & 0 & 0
		\end{pmatrix}
 $$
 where $V_1, V_2, V_3, V_4$ are Euclidean spaces. 
 \item[b)] Case of $NI =2$.
		$$
		\begin{pmatrix}
			0 &\bR& V & S^1 \\
			0 & 0 & V & S^1 \\
			0 & 0 & 0 & S^0 \\
			0 & 0 & 0 & 0
		\end{pmatrix}, 
		\begin{pmatrix}
			0 & V & W & S^1 \\
			0 & 0 & 0 & S^0 \\
			0 & 0 & 0 & 0 \\
			0 & 0 & 0 & 0
		\end{pmatrix}, 
		\ \ \ 
		\begin{pmatrix}
			0 & S^0 & W & S^1 \\
			0 & 0 & 0 & V \\
			0 & 0 & 0 & 0 \\
			0 & 0 & 0 & 0
		\end{pmatrix}, \ \ \ 
		\begin{pmatrix}
			0 & 0 & S^0_2 & S^1_2 \\
			0 & 0 & S^0_1 & S^1_1 \\
			0 & 0 & 0 & V \\
			0 & 0 & 0 & 0
		\end{pmatrix},
		$$
		where $S^0\oplus S^1$, $S^0_1\oplus S^1_1$, $S^0_2\oplus S^1_2$ are $\text{Cl}(V)$-modules;
		$$
		\begin{pmatrix}
			0 & 0 & V_2 & S^1_2 \\
			0 & 0 & V_1 & S^1_1 \\
			0 & 0 & 0 & S^0_1\simeq S^0_2 \\
			0 & 0 & 0 & 0
		\end{pmatrix}, \ \ \ 
		\begin{pmatrix}
			0 & V_1 & S_2^0 & S_1^1\simeq S_2^1 \\
			0 & 0 & 0 &S_1^0 \\
			0 & 0 & 0 & V_2 \\
			0 & 0 & 0 & 0
		\end{pmatrix} \ \ \
		\begin{pmatrix}
			0 & V_1 & V_2 & S_1^1\simeq S_2^1 \\
			0 & 0 & 0 &S_1^0 \\
			0 & 0 & 0 & S_2^0 \\
			0 & 0 & 0 & 0
		\end{pmatrix}
		$$
		where $S_1^0\oplus S_1^1$ is a $\bZ/2\bZ$-bigraded  $\text{Cl}(V_1)$ and $S_2^1\oplus S_2^1$ is a  $\bZ/2\bZ$-bigraded $\text{Cl}(V_2)$-module;
		$$
		\begin{pmatrix}
			0 & 0 & V & S^1 \\
			0 & 0 & T^0 &T^1 \\
			0 & 0 & 0 & S^0 \\
			0 & 0 & 0 & 0
		\end{pmatrix},
		$$
		where $S^0\oplus S^1$ is a  $\bZ/2\bZ$-bigraded $\text{Cl}(V)$-module and $T^0\oplus T^1$ is a  $\bZ/2\bZ$-bigraded $\text{Cl}(S^0)$-module.
	
 \end{itemize}
		
	\end{theorem}
 
 Note that all multiplications in an algebra of nilpotency index 1 are trivial. However, indicated Nil-algebras of nilpotency index 1 are not isomorphic since they have different gradations. The corresponding convex homogeneous cones are also not isomorphic.

	\begin{theorem}\label{index3}
		Any Clifford Nil-algebra of rank 4 and nilpotency index 3 admits, up to the ani-transposition, one of the following forms:
		$$
		\begin{pmatrix}
			0 & V^{10} & S^{10} & S^{11} \\
			0 & 0 & S^{00} & S^{01} \\
			0 & 0 & 0 & V^{01} \\
			0 & 0 & 0 & 0
		\end{pmatrix},
		$$
		where $V^{10}$ and $V^{01}$ are Euclidean spaces, $S^{00}\oplus S^{10}\oplus S^{01}\oplus S^{11}$ is a bigraded $\Cl(V^{10}\otimes V^{01})$-module; 
		$$	
		\begin{pmatrix}
			0 & \bR & V & S^{1} \\
			0 & 0 & V & S^1 \\
			0 & 0 & 0 & S^0 \\
			0 & 0 & 0 & 0
		\end{pmatrix},
		$$
		where $S^0\oplus S^1$ is a $\Cl(V)$-module and bilinear maps $\bR\times V \to V, \ \bR\times S^1 \to S^1 $ coincide with standard multiplications;
		$$
		\begin{pmatrix}
			0 & \C & \C & \C^n \\
			0 & 0 & \C & \C^n \\
			0 & 0 & 0 & \C^n \\
			0 & 0 & 0 & 0
		\end{pmatrix}, \ \ \ \ \ 
		\begin{pmatrix}
			0 & W & \bH & \bH^n \\
			0 & 0 & \bH & \bH^n \\
			0 & 0 & 0 & \bH^n \\
			0 & 0 & 0 & 0
		\end{pmatrix}, 
		$$
		where $W\subset \bH$ is a vector subspace.

	\end{theorem}
	 \section{Applications of homogeneous convex cones}
	\subsection{Applications to differential geometry and
supergravity}
The homogeneous convex cones are a source of construction of 
(affine and projective) special real, K\"ahler, hyperk\"ahler, and quaternionic K\"ahler manifolds that serve as the target spaces of the scalar multiplet for N=2 supergravity in dimensions D=5, 4, and 3 respectively.

\subsubsection{Homogeneous convex cones and Riemannian homogeneous geometry}
Let $C \subset V =\bR^n$ be a homogeneous convex cone and $G_o \subset Aut(C) $ a subgroup of automorphisms with a codimension one orbit $H_h = G_o y_o \subset C$ 
 that is the level set $H=H_h = \{ h(y) =1,\, y \in C\}$ of a homogeneous degree $d$ smooth function $h(y)$. Assume that $ g_H = -Hess\log h|_{TH} = - Hess h|_{TH}$ is a Riemannian metric on $H$.
 The metric $g_H$ is naturally extended to a homogeneous Lorentz metric $g$ in the cone.
 We call such a hypersurface a {\bf homogeneous Hessian hypersurface}. If, moreover, $h(y) $ is a cubic polynomial, the hypersurface $H_h$ is called a {\bf homogeneous special real manifold} in the sense of de Wit and van Proeyen. Note that the characteristic function $\chi$ of a homogeneous convex cone defines a homogeneous Hessian hypersurface and the determinant function $d(y)$ of a dual rank 3 Vinberg cone defines a homogeneous special real manifold, see \cite{alekseevsky2021special}, \cite{alekseevsky2021special1}.
 
 A homogeneous Hessian hypersurface $(H_h = G_o y_o) \subset C$ in the cone $C = \bR^{>0} H_h$ defines a complex domain
 $$D= \bR^n + i C = \{z = x+iy, x \in \bR^n, y \in C\} \subset \bC^n$$
 with a transitive group of holomorphic transformations $\bR^n \times G$ and the K\"ahler metric 
 $$
 g_D=\frac{\partial^2 k(z, \bar{z})}{\partial z^i\partial \bar z^j} dz^i d \bar{z}^j
 $$
 defined by the K\"ahler potential $k(z,\bar{z}) = -8h(\text{Im} Z) $
see \cite{cortes2017class}.

The K\"ahler manifold $(D,g_D)$ is called a {\bf Siegel domain (of the first kind)}. It is a generalization of the Poincar\'e half-plane model of the Lobachevski plane and admits a holomorphic realization as a bounded holomorphic domain. This {\bf holomorphic realization of Siegel domain}is defined as the projectivisation
$$\mathbb{P}\hat{D} = \hat{D}/\bC^* \subset \mathbb{CP}^{n+1} $$
 of the conical domain
$$ \hat{D} = \{ z^0 (1, z) , z^0 \in \bC^*,\, z \in D \}$$
 in the space $\bC^{n+1 } = \bC \oplus \bC^n$ with coordinates $(z^0, z^1, \cdots,z^n)$.
\subsubsection{r-map and c-map, special real, K\"ahler and quaternionic K\"ahler manifolds} 
 The physicists introduced {\bf special real manifold } as the target space of the scalar multiplet in
$N=2$ supergravity in dimension $D=5$. Similarly, they found that the target space of the scalar multiplet of $N=2$ supergravity in
dimension $D=4$ (resp., $D=3)$) are described by (affine or projective)
{\bf special K\"ahler manifold} (respectively hyper-K\"ahler or quaternionic K\"ahler manifolds with negative scalar curvature) 

Moreover, the dimensional reduction of supergravity from dimension $D=5$ to $D=4$ induces a
map $r : H_h \to r(H_h)$ \cite{de1992special} from a $D=5$ target space $H_h$ to the 
corresponding $D=4$ target space $r(H_h)$ , that is a special K\"ahler manifold. There are two versions of the r-map: rigid and supergravity. The image of the rigid (resp. supergravity) r-map is an affine (resp. projective) special K\"ahler manifold $K=r(H_h$ (resp., $\bar{K} =\bar{r}(H_h)$).
Similarly, the dimensional reduction from $D=4$ to $D=3$ induces a c-map , rigid ( resp., supergravity). that associates a hyperk\"ahler manifold $Q= c(K)$ 
( resp., a quaternionic K\"ahler manifold $\bar{Q}= \bar{K}$ with negative scalar curvature) to an 
affine special K\"ahler manifold $K$ (resp., projective special K\"ahler manifold $\bar{K}$.
The supergravity r and c maps preserve the homogeneity, i.e. they map homogeneous manifolds into homogeneous manifolds. 

These constructions have a purely differential geometric description, which we briefly discuss for homogeneous manifolds. 

 Let $H_h \subset C$ be a special homogeneous real manifold, defined by a homogeneous cubic $h$ 
 in a homogeneous convex cone $C$. The rigid r-map $r : H_h \to r(H_h)$ converts hypersurface $H_h$ into an affine homogeneous special 
 K\"ahler manifold $r(H_h ) = (M,g, J, \n)$ , that is a K\"ahler manifold $(M,g,J)$ with a flat torsion-free connection $\n$, that preserves the K\"ahler form $\omega = g \circ J$ and satisfies 
 the condition $(\n_X J) = (\n_Y J)X$. It is shown in \cite{alekseevsky2009geometric} that the 
 manifold $ H_h$ is identified with the tangent bundle $TH_h$ with the naturally defined special K\"ahler structure.
 
 The supergravity r-map $\hat{r} ; H_h \to \hat{r}(H_h) = \hat{K}$ is the special case of the projective realization of the Siegel domain, described above.
 
 The rigid c-map converts a $2n$ dimensional affine special (pseudo)K\"ahler manifold $(K, g_K, J,\n)$ into a (pseudo)hyperK\"ahler manifold $Q= c(K)$. The manifold $Q$ can be identified with the cotangent bundle $Q = T^*K,$ which induces hyperK\"ahler structure $( g_Q, J_1, J_2, J_3 = J_1J_2$. The decomposition 
 $$ T^*K = (T^*)^vK \oplus (T^*)^h(K)
 $$
 of the bundle $T^*K$ into vertical and horizontal subbundles, defined by the flat connection $\n$, allows us to extend the complex structure $J$ and the K\"ahler form $\omega = g_K\circ J$ into anticommuting complex structures $J_1, J_2$ and the metric $g_K$ into a hyperK\"ahler metric $g_Q$.
 
 The construction of a supergravity c-map $\bar{c}z: \bar{K} \to \bar{c}(\bar{K}) = \bar{Q}$, which converts a projective special K\"ahler manifold $\bar{K}$ into a quaternionic K\"ahler manifold $\bar{Q}$, is more complicated.
 In \cite{alekseevsky2013conification} a purely differential geometric construction of this c-map was proposed.
 It is based on a procedure of conification, which converts a (pseudo)hyperK\"ahler $4n$-dimensional manifold $(Q,g_Q, J_1,J_2,J_3)$ with an appropriate Killing vector field $\xi$ into a conical (pseudo)hyperK\"ahler manifold $\hat{Q}$ of dimension $4n+4$, which is the (total space of ) the Swan bundle 
 $\pi : \hat{Q} \to \bar{Q} = P\hat{Q} := \hat{Q}/(\bH^*/\bZ_2)$. The quaternionic projectivization of $\hat{Q}$ is a (pseudo)quaternionic K\"ahler manifold $\bar{Q}$.
 
 The supergravity c-map $\bar{c} : \bar{K} \to \bar{Q} = P \hat{Q}$ is defined as the 
 quaternionic K\"ahler manifold, which is the projectivization of the conical hyperK\"ahler manifold ${\bar{c}(\bar{K})}$, obtained by the conification of a (pseudo)hyperK\"ahler manifold $\bar{c}(\bar{K})$ of signature $(2,2n)$, associated with the affine c-map.

\subsubsection{ Rank 3 Vinberg cones and homogeneous special real K\"ahler and quaternionic K\"ahler manifolds}
 The classification of homogeneous special real manifolds $H_h$ was done by de Wit and Van Proeyen \cite{de1992special} and by another method, by V.Cortes \cite{cortes1996homogeneous}. 
In \cite{alekseevskii1975structure}
 D.V. Alekseevsky classified quaternionic K\"ahler manifolds with a transitive solvable group of isometries. 
 They bijectively corresponds to $\bZ_2$ graded Clifford modules and have rank $\leq 4$. See also \cite{cortes1996homogeneous} for some corrections.
 Alekseevsky conjectured that such manifolds exhaust all homogeneous quaternionic K\"ahler manifolds of negative scalar curvature, and, more generally, any homogeneous Einstein manifold of negative scalar curvature is a solvmanifold, i.e., it admits a transitive solvable Lie group of isometries. Both conjectures have been recently proved by Ch. B\"ohm and R. Lafuente \cite{bohm2022homogeneous}
S. Cecotti \cite{cecotti1989geometry} proved that the associated Alekseevsky quaternionic K\"ahler spaces, different from symmetric quaternionic K\"ahler manifolds
$$
\bH P^m = Sp(1, m)/Sp(1) \times Sp(m) \ \ \ \text{and} \ \ \ Gr_{m+2}(\bC^{m+2})=U(m, 2)/S(U(m) \times U(2)) 
$$
are special, i.e. are in the image of the supergravity c-map.\\
Recently in \cite{matrix} it was proven that all rank 4 Alekseevsky spaces bijectively corresponds to rank 3 special Vinberg cones, that is, they are in the image of the $q = c \circ r$ map.\\
In the paper \cite{alekseevsky2021special1}, the authors use the rank 3 Vinberg cones to determine an explicit formula for the Bekenstein-Hawking entropy of a static, spherically symmetric and asymptotically flat BPS extremal black hole in terms of its electric and magnetic charges in the ungauged N = 2 D = 4 supergravity theory, under the assumption that the scalar multiplet ( more precisely, the scalar manifold of the vector multiplets) is homogeneous.
 
Due to the Bekenstein-Hawking entropy-area formula, Ferrara,
Kallosh and Strominger’s BPS algebraic attractor equations, and the 
 M. Shmakova \cite{shmakova1997calabi}
 formula, the problem reduces to inversion of the quadratic 
diffeomorphic map
$$ d: C'(S) \to C'(S)^*; X \to d(X, X.\cdot ) \in C(S) $$
from the dual rank 3 Vinberg cone $C'(S)$ to the adjoint cone $C'(S)^*$
where $d(X,Y,Z)$ is the polarization of the cubic determinant 
$d(X) = d(AA^*) = (a_{11} a_{22} a_{33})^2 $.

\subsubsection{Hyperbolic affine spheres and Monge-Ampere equation} 
 Let $S \subset \bR^n$ be a hypersurface transversal to a radiant vector field $\rho=\sum x^i \partial_{\dot{ x}^i}$. Then for $x\in S$ the tangent space $T_x \bR^n = \bR^n$ is decomposed as 
 $T_x \bR^n = \bR \rho \oplus T_xS$.
 The second fundamental form $g(X,Y)$ of $S$ is defined by $$ [\tilde{X}, \tilde{Y}] = g(X,Y) \rho + [\tilde{X}. \tilde{Y}]_{TS}$$
 where $ \tilde{X},\tilde{Y} $ are extensions of tangent vector fields $X,Y \in \mathcal{X}(S)$ to a neighborhood of $S$.
 If $g$ is positive definite, it is called the Blashke metric. If, moreover,
 the Riemannian volume form $\text{vol}^g$ coincides with the volume form 
 $$(\iota_{\rho}) dx|_S = \iota_{\rho} (dx^1 \wedge \dots \wedge dx^n)|_S , $$ 
 then $(S,\rho)$ is called a {hyperbolic affine hypersphere}. Affine hyperspheres are introduced by W. Blashke (\cite{blaschke1923vorlesungen}). See \cite{loftin2008survey} for a survey on affine spheres.\\
 
 Any complete hyperbolic affine hypersphere is the level hypersurfaces $u(x)=1$ of a solution $u(x)$ of a (unique global) Monge-Ampère equation on a convex cone $C$ (\cite{calabi1972complete}, \cite{cheng1986complete}). Conversely, any convex cone determines an affine hyperbolic sphere (\cite{gigena1981conjecture}.\\
 For a homogeneous convex cone $C$ , the solution $u(x)$ is the characteristic function up to a constant multiplier (\cite{sasaki1980hyperbolic}).

 \subsection{Applications to Hessian geometry, statistics and convex programming}

\subsubsection{ Homogeneous convex domains and left symmetric algebras } 
Homogeneous convex cone is a special case of homogeneous (w.r.t. an affine group ) convex domains in an affine space $\bA^n = \bR^n$ . Any bounded homogeneous convex domain is obtained as an intersection of some homogeneous convex cone with an affine hyperplane. The investigation of homogeneous convex domains was initiated by Koszul \cite{koszul1965varietes}. The classification of homogeneous convex domains reduces to an algebraic problem of description of some class of left symmetric algebras (also called Koszul-Vinberg algebras and pre-Lie algebras) defined by Cayley, Koszul, and Vinberg. Left symmetric algebras have many applications. See \cite{burde2006left}, \cite{manchon2011short}, and \cite{bai2021introduction} for a survey on LSAs.

An algebra $(\mathcal{G},\cdot)$ is called a {\bf left symmetric algebra (LSA)} if the multiplication satisfies the identity
$$
(x\cdot y)\cdot z - x\cdot (y\cdot z)=(z\cdot x)\cdot y - z\cdot (x\cdot y)
$$
for all $x,y,z\in\mathfrak{g}$. The commutator 
$$
[x,y]=x\cdot y-y\cdot x
$$
defines in $\mathcal{G}$ the structure of a Lie algebra $\mathfrak{g}$. If $G$ is a Lie group with 
$\text{Lie}(G) = \gg$ , then the formula $\n_x y=x\cdot y$ defines in $G$ a left-invariant torsion-free flat connection, where $x,y \in \gg$ is considered as left invariant vector fields. The inverse claim is also true.

According to Vinberg \cite{vinberg1963theory}, the classification of homogeneous convex domains without straight lines reduces to a classification of {\bf compact left symmetric algebras (clans)} $\mathcal{G}$, 
such that the adjoint Lie algebra $\ad_{\gg}$ is triangular and there is a linear form $s \in \mathcal{G}$ such that $g(x,y) := s(x \cdot y)$ is a Euclidean metric.
 In particular, homogeneous convex cones corresponds to clans with unit $e \in \mathcal{G},\, e\cdot x = x \cdot e$.
 The theory of T-algebras, that provides a description of clans with unit, reduces the classification of homogeneous convex domains to the description of sections of homogeneous convex cones $C \subset \text{Herm}_n $ by affine hyperplanes, which are homogeneous.\\
 In particular, our classification of rank 4 Vinberg cones provides a large new class of compact left-symmetric algebras.
 
 \subsubsection{Homogeneous Hessian manifolds}

 A flat affine manifold is a manifold $M$ endowed with a flat torsion-free connection $\nabla$. A Hessian manifold is a flat affine manifold $(M,\nabla)$ together with a Riemannian metric $g$ locally equivalent to a Hessian of a function $g=\nabla^2f$.
 	Any K\"ahler metric is locally expressed as a complex Hessian of a function. Thus, Hessian geometry can be considered as a real analogue of K\"ahler geometry. Moreover, if $M$ is a Hessian manifold, then $TM$ admits a K\"ahler structure. In the case of a convex cone $C$, the tangent bundle $TC\simeq \bR^n\oplus \sqrt{-1}C $ is biholomorphic to a bounded domain in $\C^n$. The K\"ahler metric on $TC$, constructed by the canonical Hessian metric on $C$, coincides with the Bergman metric (see \cite{shima2007geometry}).

	Any convex cone $C \subset \bR^n$ admits a canonical invariant under linear automorphisms Hessian metric $g_C$. We can consider any convex domain $U$ as a section of a convex cone $C$. The restriction of the canonical Hessian metric $g_C$ on $U$ is an invariant under affine automorphisms Hessian metric. Thus, any homogeneous convex domain without straight lines admits a structure of a homogeneous Hessian manifold. Conversely, if there are no non-constant affine maps from $\bR$ to a homogeneous Hessian manifold $M$, then the universal covering of $M$ is an affine homogeneous convex domain (\cite{shima1980homogeneous}).
 
 Moreover, if a homogeneous Hessian manifold $M$ is symmetric or admits a transitive reductive group of affine automorphisms, then the universal covering of $M$ is a direct product of a Euclidean space and a homogeneous convex self-dual cone $C$, which is a symmetric space with the canonical metric $g_C$(\cite{shima1977symmetric},\cite{shima1980homogeneous}).
	

 \subsubsection{Information geometry and statistics}
Different families of probability distributions on convex homogeneous cones are studied in statistics. The homogeneous convex cones are especially important because many important invariants (e.g. the characteristic function) may be easily calculated explicitly.

 In particular, the characteristic function $ \chi(\xi) =\int_{C^*} e^{-\langle x,\xi \rangle } dx$ of a convex cone $C$ defines the family of probability distributions $\mu_\xi = \frac{e^{-\langle x,\xi \rangle }} {\chi(\xi)}$, parametrized by points of the dual cone $C^*$. 
 
 It is called an exponential family, see \cite{morris2014natural} .
The canonical Hessian metric $g_C = \partial^2 \log \chi(\xi)$, invariant w.r.t the automorphism group, is defined for any convex cone $C$. It coincides with the Fisher-Rao metric of the exponential family.
 The exponential families play a distinguished role in statistics and information geometry. In particular,the exponential distributions can be characterized by the principle of maximum entropy. Any exponential distribution is the distribution with maximal entropy among all distributions with fixed mean value. 
 
The generalization of the exponential family is the Wishart family, defined by Wishart for the cone of positive definite real matrices in \cite{wishart1928generalised} and generalized by Andersen to any homogeneous convex cone \cite{andersson2004wishart}.

The Riesz family of probability distributions, introduced in \cite{riesz1949integrale} for the Lorentz cone in $\bR^{1,n}$, was generalized to an arbitrary homogeneous convex cone in \cite{gindikin1964analysis}. For the study of Riesz and Wishart families on homogeneous convex cones, see also \cite{ishi2000positive} and \cite{hassairi2005beta}. 

The classical Gaussian distributions on the space of positive symmetric definite matrices are studied, for example, in \cite{said2017riemannian}. The generalization of these distributions to
 Riemannian symmetric spaces of non-positive curvature (in particular, to symmetric homogeneous convex cones) are studied in \cite{said2022gaussian}. 
 For applications of homogeneous convex cones in statistics, see the book \cite{schneider2022convex}.

\subsubsection{Homogeneous convex cones and convex programming}
 Homogeneous convex cones and the characteristic function $\chi(x)$ are widely used in convex programming. A. Nemirovski and Yu. Nesterov \cite{nemirovski2006advances} developed a theory of conic optimization (based on the internal-point method) and proposed an effective method to find the minimum of a convex function $f(x)$ in a convex cone $C$. The main idea is to use the characteristic function $\chi$ as a barrier function (tending to infinity on the boundary and preventing the iterative process from approaching the boundary). It is applicable primarily to homogeneous cones, for which this function is explicitly described. 
The approach consists of replacing the convex function $f(x)$ in a convex cone by the family
$F_t(x) = t f(x) +\chi(x) $. 
The path $x_t = \mathrm{argmin}{F_t}(x)$ of the minimizers
converges to the minimizer $x_o =\mathrm{argmin}f(x)$ for $t \to \infty$ of the initial function. 
 The Newton method ( or other iterative procedures) allows one to approximate this path by a sequence $x(i)$ of points that approximate the minimizer $x_{t_i}$ for a chosen increasing sequence $(t_i) $ with $ t_{i+1}=\mu t_i$ such that for sufficiently big $t\_i$ the point $x(i)$ gives a good approximation for the initial minimizer $x_0$. 
The sequence is initialized by a certain nontrivial two-phase initialization procedure that uses the cone $C$ and the dual cone $C^*$, see \cite{nemirovski2006advances}.

There is a large literature devoted to the development and generalization of this approach, in particular, to convex domains, see \cite{nesterov1994interior}, \cite{ben2001lectures}.
	
 \section{ Nil-algebras and homogeneous convex cones}
	
	\subsection{Isometric maps}
	Let $X,Y, Z$ be Euclidean vector spaces. 
	\begin{Def}
		A bilinear map $$ \mu : X \times Y \to Z,\, \ \ \ (x,y)\to \mu(x,y)= x\cdot y $$
		is called an {\bfseries isometric map} if
		the norm $|x\cdot y|= |x| \cdot |y|.$
	\end{Def} 

 Hurwitz posed the problem of finding all isometric maps in 1898 (\cite{hurwitz1898ober}. This problem is still open in the general case (see \cite{shapiro2011compositions}).

 \begin{Example}\label{EX2}
		$X = Y= Z=\bA$, where $\bA$ is a normed division algebra. Then the multiplication in $\bA$ is an isometric map.
	\end{Example}

	\begin{Example}\label{EX}
		Let $X=V$ be a Euclidean space, $\Cl(V)$ the corresponding Clifford algebra, ${Y=S^0}$, $Z=S^1$ (or $Y=V,X=S^0,Z=S^1)$, where $S=S^0\oplus S^1$ is a $\bZ_2$-graded $\Cl(V)$-module equipped with a
 $Spin(V)$-invariant Euclidean metric $g_S$ with $S^0, \perp S^1$, 
 called an {\bf admissible metric.}
 Then the Clifford multiplications $V \times S^0 \to S^1$, $(v,s^0)\mapsto vs^0$ (or $S^0 \times V \to S^1, (s^0,v) \mapsto vs^0$) 
 are isometric maps.
 
 The theory of Clifford modules was founded in \cite{atiyah1964clifford}, see also \cite{lawson2016spin}.
 The description of an admissible (pseudo)Euclidean metric $g_S$ for $\bZ_2$ graded Clifford $Cl(V)$ modules for (pseudo)Euclidean vector spaces $V$ was given in \cite{alekseevsky1997classification}.\\
	\end{Example}
{\bf We will always assume that a Clifford module $S = S^0 + S^1$ is $\bZ_2$ graded and is endowed with an admissible Euclidean metric.}

	\begin{theorem}[ \cite{hurwitz1922komposition}] \label{Hurwitz}
			Let $\mu: X\times Y \to Z$ be an isometric map and $\dim X=\dim Y =\dim Z$. Then there is an identification $X\simeq Y\simeq Z\simeq \bA$, where $\bA\in \{\mathbb{R},\C,\mathbb{H},\mathbb{O}\}$, such that $\mu$ coincides with the multiplication in $\bA$. 
	\end{theorem}
	
	\begin{theorem}[\cite{lee1948theoreme},\cite{dubisch1946composition}]\label{ABSh}
		Let $\mu: X\times Y \to Z$ be an isometric map and $\dim Y =\dim Z$ (or $\dim X=\dim Z$). Then $X,Y,Z$ are as in Example $\ref{EX}$.\\
	\end{theorem}
	
	\begin{rem}
		We have $\dim V =\dim S^0=\dim S^1$ for an indecomposable module $S^0\oplus S^1$ if and only if $\dim V=1,2,4,8$ and $S$ is an indecomposable $\Cl(V)$-module (\cite{atiyah1964clifford}). In this case, according to Theorem \ref{Hurwitz}, the ordered triples $(V,S^0,S^1)$, $(S^0,V,S^1)$, and $(\bA,\bA,\bA)$ (together with the corresponding isometric maps) are isomorphic. 
	\end{rem}
	 
	\begin{Lem}\label{lemmaXYZ}
		Let $k\in \{1,2,4\}$, $\mu: X\times Y \to Z$ be an isometric map with $\dim X \leq\dim Y =\dim Z=k$, $\bA$ a division algebra with $\dim \bA = k$, $W\subset \bA$ be a vector subspace with $\dim W=\dim X$. 
		\begin{itemize}
			\item[a)] The triple $X,Y,Z$ with isometric map $\mu$ is isomorphic to $(W,\bA,\bA)$ with the standard multiplication map $W\times \bA \to \bA$.
			\item[b)] Suppose that $1\in W\subset \bA$. 	Then for any $e_X\in X, e_Y\in Y$ there exist isometries 
			$$
			\varphi_X : X \to W, \ \ \ \varphi_Y : Y \to \bA, \ \ \ \varphi_Z : Z \to \bA
			$$
			such that 
			$$
			\varphi_X(e_X)=1, \ \ \ \ \varphi_Y(e_Y)=1, \ \ \ \ \varphi_Z \left(\mu (x,y) \right) =\varphi_X (x)\varphi_Y (y).
			$$ 
		\end{itemize}

	\end{Lem}
	
	\begin{proof}
		a) Suppose $\dim X\ \le 3$. Let $V$ be a Euclidean space with $\dim X=\dim V \le 3$. Then there exists a unique indecomposable $\bZ/2\bZ$-graded $\Cl(V)$-module (\cite{atiyah1964clifford}) $S^0\oplus S^1$. According to Theorem \ref{ABSh}, the triples $(X,Y,Z)$ and $(W,\bA,\bA)$ (together with the corresponding isometric maps) are both isomorphic to $(V,S^0,S^1)$. Thus, $(X,Y,Z)$ and $(W,\bA,\bA)$ are isomorphic.

 Suppose $\dim X=4$. Hence, we have $\dim X=\dim Y=\dim Z$. According to Theorem \ref{Hurwitz}, the triple $(X,Y,Z)$ is isomorphic to $(\bH,\bH,\bH)$.
		
		b) According to item a) there exists an isomorphic map $\psi_X,\psi_Y,\psi_Z$ between $(X,Y,Z)$ and $(W,\bA,\bA)$
			such that 
			$$
			\psi_Z \left(\mu (X,Y) \right) =\psi_X (x)\psi_Y (y).
			$$ 
 Then the maps
		$$
		\varphi_X(x)=\left(\psi_X(e_X)\right)^{-1} \psi_X (x), \ \ \ \ \ \varphi_Y(y)=\psi_Y(y) \left(\psi_X(e_Y)\right)^{-1},
		$$
		$$ 
		\varphi_Z(z)=\left(\psi_X(e_X)\right)^{-1} \psi_Z(z) \left(\psi_X(e_Y)\right)^{-1}
		$$
		satisfy the conditions of item b).	\end{proof}

	\subsection{(Quasi)Nil-algebras}

 \subsubsection{Definitions}
	
	Let $\cn = \sum_{i<j =1}^n \cn_{ij}$ be an orthogonal decomposition of a Euclidean vector space $(\cn, g = \langle ., . \rangle)$.
	We call $\mathcal{N}$ a bigraded Euclidean space.
	It is convenient to
represent an element $A = \sum_{i,j =1}^m a_{ij}$ as $m \times m$ matrices $A = ||a_{ij}||$. A system of bilinear maps
	\be \label{matrix multiplication}
	\cn_{ij} \times \cn_{jk} \to \cn_{ik} 
	\ee
	 defines in $\cn$ the structure of an algebra (called {\bf generalized nilpotent matrix algebra}) with the standard rule of matrix multiplication.. 
 
	 The {\bfseries nilpotency index} of a nilpotent associative algebra is the number $NI(\cN)=d$, where $d-1$ is the maximal number of elements $A_1, \cdots A_{d-1} \in \cN$ with non-zero product
 $A_1 \cdot A_2 \cdot \cdots \cdot A_{d-1} \neq 0$.
		
 \bd
		i) A {\bfseries (quasi)Nil-algebra} of rank $n$ is a bigraded Euclidean space $\cn = \sum_{1\le i<j\le n} \cn_{ij}$,
 with a matrix multiplication
 $$\cn \times \cn \to \cn ,\,(A, B) \mapsto C=AB$$
			where $ C= c_{ik} = \sum_j a_{ij}b_{jk}$
			defined by a system of isometric maps
			$$ \mu_{ijk}: \cn_{ij} \times \cn_{jk} \to \cn_{ik},\,\, i<j<k. $$
 ii) \cite{vinberg1963theory} \label{Vinberg} 	 A {\bf Nil-algebra} is an associative (quasi)Nil-algebra satisfying the Vinberg condition: 
			for any $a_{ik} {\in} \cn_{ik}$, $a_{jk} {\in} \cn_{jk}$ with $i < j$ 
	$$
 \text{if} \ \ \ \ \langle a_{ik}, \cn \cdot a_{j k}\rangle {=}\{0\} \ \ \ \
 \text{then} \ \ \ \
 		\langle \cn \cdot a_{ik}, \cn \cdot a_{jk} \rangle{=} \{0\}\ .
 $$
 \ed 
 \begin{Remark}\label{rank4vinbercondition}
i) If $n= 4$ then the Vinberg condition reduces to the following condition:
for any $a_{24} {\in} \cn_{24}$, $a_{34} {\in} \cn_{34}$ 
			\begin{equation} 
				\label{Vinberg condition}
 \text{if} \ \ \ \
 \langle a_{24}, \cn_{23} \cdot \a_{3 4}\rangle {=}\{0\} \ \ \text{then} \ \
				\langle \cn_{12} \cdot a_{24}, \cn_{34} \cdot a_{34} \rangle{=} \{0\}\ .
			 \end{equation}
 ii) Suppose $n= 4$ and the nilpotency index $NI(\cN)$ of a quasiNil-algebra $\cN$ is $\leq 2 $, then the algebra $\cN$ is associative and there exists $1\le i< j\le k$ such that $\cN_{ij}= \{0\}$. Therefore, Condition \eqref{Vinberg condition} is satisfied and $\cN$ is a Nil-algebra. 
 \end{Remark}

 \subsubsection{ Isomorphism and equivalence}
\bd
 Two quasiNil-algebras $\cN$ and $\cN'$ are {\bf isomorphic} if there exists an isomorphism of algebras $\varphi: \cN \to \cN'$ (that preserves the bigradation and the metric).
\ed
 We will consider the ( standard ) matrix 
 realization of the (quasi-)Nil-algebra $\cn$ as a nilpotent associative algebra of upper 
 triangular matrices 
 of the form
	\begin{equation}
		A = \begin{pmatrix}
			0 & a_{12} & a_{13} & \ldots & \ldots & a_{1\, n-1} & a_{1\,n}\\
			0 & 0 & a_{23} & \ldots & & a_{2\, n-1} & a_{2\,n}\\
			0 & 0 & 0 & \ldots & & a_{2\, n-1} & a_{3\,n}\\
			\vdots & \vdots & \ddots & \ddots & \ddots &\vdots & \vdots\\
			0 & 0& 0 & \ldots & \ldots &0 & a_{n-1\,n}\\
			0 & 0& 0 & \ldots &\ldots & 0 & 0
		\end{pmatrix}
	\end{equation}
	with matrix multiplication, defined by the isometric maps $\cn_{ij} \times \cn_{jk} \to \cn_{ik}$.
 The canonical metric $g$ is given by $g(A,B) = \tr AB^t = \sum \langle a_{ij}, b_{ij}\rangle $.
 \begin{Def}
 \label{admissible}
	i) A permutation $ \sigma \in S_n$ that changes the bigradation $\cn =\sum \cn_{ij}$ 
 of a quasiNil-algebra into $\cn =\sum \cn_{\s(i) \s(j)}$ is called {\bf admissible}, if the algebra with the new bigradation $\cN_{\s(i) \s(j)}$ , called a {\bf renumbering} of $\cN$, is again a quasiNil-algebra .\\
 ii) Two quasiNil-algebras are called {\bf equivalent}, if they become isomorphic after a renumbering.\\
 iii) The {\bf dual } quasiNil-algebra $\cn^{\tau}$ as a quasiNil-algebra $\cN$ 
 obtained from $\cn$ by changing the gradation to $N'_{ij} := N_{n+1-j \,n+1 -i}$.
	In matrix notation, the dual algebra is obtained from $\cn$ by the anti-transposition $\tau$, that is, the reflection across the anti-diagonal. \\
 iv) A quasiNil-algebra $\cn$ is called {\bf Clifford}, if all non-zero isometric maps 
		$$
		\cn_{ij} \times \cn_{jk} \to \cn_{ik}
		$$
 are Clifford, i.e. are identified with a Clifford multiplication
 $V \times S^0 \to S^1$
 of a graded Clifford $Cl(V)$ module $S = S^0 \oplus S^1$, where
 $V=\cn_{ij},\,S= \cn_{jk} \oplus \cn_{ik} $
 or $S= \cn_{ij} \oplus \cn_{ik},\, V=\cn_{jk},$.\\
 Equivalently, for any non-trivial map
		$\cn_{ij} \times \cn_{jk} \to \cn_{ik}
		$ we have $\dim \cn_{ij}=\dim \cn_{ik}$ or $\dim \cn_{jk}=\dim \cn_{ik}$.
	\end{Def}
		\subsubsection{Indecomposable quasiNil-algebras}
 \begin{Def}
 i) The {\bf direct sum} ${\cN}= \cN'\oplus \cN''$ of rank $n'$ quasiNil-algebra $\cN'$ and a rank $n''$ 
 graded quasiNil-algebra $\cN''$ is a quasiNil-algebra with the bigradation \
 \begin{align*}
 \cN_{ij}=\cN'_{ij} \ \ \ &\text{for} \ \ \ 1\le i<j <n' \\
 \cN_{n'+i \ n'+j}=\cN''_{ij} \ \ \ &\text{for} \ \ \ 1\le i<j <n''
 \end{align*}
 ii) A (quasi)Nil-algebra $\cN$ is {\bf decomposable} if it is equivalent to a direct sum of two (quasi)Nil-algebras .
 \end{Def}
 In other words, a (quasi)Nil-algebra $\cN$ is decomposable if it is equivalent to a (quasi)Nil-algebra of block diagonal matrices.
A Nil-algebra is indecomposable if and only if the corresponding cone is indecomposable, i.e. is not linearly isomorphic to a direct product of two cones.
The classification of (quasi)Nil-algebras reduces to the classification of indecomposable ones.
		\smallskip
	\subsection{Reconstruction of Vinberg T-algebra and a homogeneous convex cone from a Nil-algebra}
	\subsubsection{Solv-algebras and Vinberg groups}
	
	Given a Nil-algebra $\cN$ of rank $n$ over $\bK$, the matrix multiplication between elements in $\cN$ naturally extends to the matrix multiplication between upper triangular matrices of the form
	\begin{equation*} \label{V''}
		\begin{split}
			& A = D + N
			= \left(\smallmatrix
			a_{11} & a_{12} & a_{13} & \dots & \ldots & a_{1 n-1} & a_{1n}\\
			0 & a_{22} & a_{23} & \ldots & & a_{2 n-1} & a_{2n}\\
			0 & 0 & a_{33} & \dots & & a_{2 n-1} & a_{3n}\\
			\dots & \dots & \dots & \dots & \dots &\dots & \dots\\
			0 & 0& 0 & \dots & \dots & a_{n-1 n-1} & a_{n-1n}\\
			0 & 0& 0 & \dots &\dots & 0 & a_{nn}
			\endsmallmatrix\right)\\[10pt]
			& \text{with} \ N \in \cn\ ,\ D = \operatorname{diag}(a_{11}, \dots, a_{nn}) \ \text{with}\ a_{ii} {\in} \bR\ ,
		\end{split}
	\end{equation*}
	by means of the natural maps
	$$
	( a_{ii}, a_{ik})\mapsto a_{ii} a_{ik}.
	$$
	The associative algebra $\cg {=} \cg(\cn) {=} D {+} \cn$ of the matrices of the form \eqref{V''}
is called {\bfseries the Solv-algebra associated with $\cn$}.
	The commutator $[A,B] = AB -BA$ defines in $\cg$ the structure of a solvable Lie algebra , which will be denoted by $\gg = (\cg, [.,.])$. 
	\par
	
	\bd 
	The {\bf Vinberg group} $G = G(\cn)= \{A = ||a_{ij}||\in \cg,\, a_{ii}>0 \}$ is the (simply connected)
	connected solvable Lie group of invertible matrices with the Lie algebra $\gg$.
 \ed
	Note that the Vinberg groups $G(\cn)$, $G(\cn^{\tau})$, associated with $\Nil$-algebras $\cn, \cn^{\tau}$ 
 are isomorphic iff the algebra $\cn^{\tau}$ is equivalent to $\cn$ or its dual algebra $\cn^{\tau}$, 
	\par
	\medskip
	\subsection{The T-algebra associated with a Nil-algebra}
	Let $\cn = \sum_{i<j}\cn_{ij}$ be a Nil-algebra of rank $n$. Following Vinberg, we extend it to the full matrix algebra $ \cT= \cT(\cn)$, called T-algebra \cite{vinberg1963theory}, as follows. As a
	vector space $\cT$ consists of $m \times m$ matrices of the form
	$$
	M = \sum_{i,j =1}^m M_{ij }= \left(\smallmatrix
	a_{1\,1} & a_{1\,2} & a_{1\,3} & \ldots & \ldots & a_{1\,n-1} & a_{1\,n}\\
	a_{2\,1} & a_{2\,2} & a_{2\,3} & \ldots & & a_{2\,n-1} & v_{2\,n}\\
	a_{3\,1} & a_{3\,2} & a_{3\,3} & \ldots & & a_{2\,n-1} & v_{3\,n}\\
	\vdots & \vdots & \ddots & \ddots & \ddots &\vdots & \vdots\\
	a_{m-1\,1 } & a_{m-1\,2 }& a_{m-1\,3} & \ldots & \ldots & a_{n-1\,n-1} & a_{n-1\,n}\\
	a_{n\,1 } & a_{n\,2 }& a_{n\,3 } & \ldots &\ldots & a_{n n-1} & a_{n\,n}
	\endsmallmatrix\right)\ ,
	$$
	whose entries $a_{ij}$ are elements of the vector spaces
	$$a_{ij}\in \left\{\begin{array}{ll} \cn_{ij}& \text{if} \ i < j\ ,\\[10pt]
		\bR& \text{if}\ i = j\ ,\\[10pt]
		\cn_{ji}^* = \Hom(\cn_{ij}, \bR) & \text{if}\ i > j\ .
	\end{array} \right. $$
	The space $\cT(\cN)$ has a canonical structure of non-associative algebra with multiplication, which extends the multiplication of the associative algebra $\cg$.
 The missing bilinear products $\cn_{ij} \times \cn_{jk}$ for $i>j$ or $j>k$
 are given by 
		
		\begin{equation*}
			\begin{split}
				& (a_{ij}\cdot a_{jk})^* = a_{jk}^* \cdot v_{ij}^*,\\
				& \langle a_{ij}^* \cdot a_{ik}, a_{jk}\rangle = \langle a_{ik}, a_{ij}\cdot a_{jk}\rangle,\\
				&\langle a_{ik}\cdot a_{jk}^*, a_{ij}\rangle = \langle a_{ik},a_{ij} \cdot a_{jk}\rangle ,
			\end{split}
		\end{equation*}
		where we denote by $v_{ij}^* = \langle v_{ij}, \cdot \rangle$ the covector associated with the vector $v_{ij}$. \par
		\smallskip
		\bd
		The algebra $\cT = \cT(\cn)$ is called the {\bf T-algebra} associated with the Nil-algebra $\cn$. Two T-algebras $ \cT, \cT'$ are called {\bf equivalent} iff the corresponding Nil-algebras are equivalent.
		\ed
	

		\par
		\medskip
		
		\subsection{The space \texorpdfstring{$\cH = \text{Herm}(\cn)$}{} of Hermitian matrices and the action of the Vinberg group \texorpdfstring{$G$}{G} on \texorpdfstring{$\cH$}{H}}
		 We denote by $\cH = \Herm(\cn) = \{ X =X^* \} \subset {\mcal T}$ the subspace of Hermitian matrices, where the conjugate to $X = ||x_{ij}||$ 
 matrix is the matrix $X^* =||x^*_{ij}||$ with the entries $(x^*)_{ji}:= x_{ij}^* :=(x_{ij})^{*}$. 
 The subspace $\cH = \Herm(\cn)$ is a commutative non-associative algebra with the Jordan multiplication
		$$ \{X,Y\} = \frac12(XY + YX).$$
		 Any upper triangular matrix
		$$A = D + N = \left(\smallmatrix
		a_{11} & a_{12} & a_{13} & \ldots & \ldots & a_{1 n-1} & a_{1n}\\
		0 & a_{22} & a_{23} & \ldots & & a_{2 n-1} & a_{2n}\\
		0 & 0 & a_{33} & \ldots & & a_{2 n-1} & a_{3n}\\
		\vdots & \vdots & \ddots & \ddots & \ddots &\vdots & \vdots\\
		0 & 0& 0 & \ldots & \ldots & a_{n-1 n-1} & a_{n-1n}\\
		0 & 0& 0 & \ldots &\ldots & 0 & a_{nn}
		\endsmallmatrix\right) $$
		is 
		 naturally extended to the {\bfseries Hermitian matrix}
		$$ H_A = N^* + D + N = \left(\smallmatrix
		a_{11} & a_{12} & a_{13} & \ldots & \ldots & a_{1 n-1} & a_{1n}\\
		a_{12}^* & a_{22} & a_{23} & \ldots & & a_{2 n-1} & a_{2n}\\
		a_{13}^* & a_{23}^* & a_{33} & \ldots & & a_{2 n-1} & a_{3n}\\
		\vdots & \vdots & \ddots & \ddots & \ddots &\vdots & \vdots\\
		a_{1 n-1}^* & a_{2 n-1}^* & a_{3 n-1}^* & \ldots & \ldots & a_{n-1\,n-1} & a_{n-1\,n}\\
		a_{1 m}^* & a_{2 n}^* & a_{3 n}^* & \ldots &\ldots & a^*_{n-1\,n} & a_{n\,n}
		\endsmallmatrix\right) \in \bT(\cn)\ .$$
		\par
		\medskip
		The next theorem collects some crucial results of Vinberg \cite{vinberg1963theory} on homogeneous cones in spaces of Hermitian matrices.
		\bt[\cite{vinberg1963theory}] \hfill\par
		\begin{itemize}
			\item The map
			\begin{equation}\label{ve}
				L: \gg \times \Herm(\cn) \longrightarrow \Herm(\cn)\ ,\qquad L_A(X)= \dot A X + X \dot A^*
			\end{equation}
			
			is a faithful linear representation of the Lie algebra $\gg = \gg(\cn)$ of the Vinberg group $G = G(\cn)$ in the vector space $\Herm(\cn)$.
			\item The Lie algebra representation $L$ integrates to a Lie group representation
			$$ \rho: G \to \GL{\cH}$$ of the Vinberg group $G$ on $\cH=\Herm(\cn)$
			\footnote{Such $G$-action is determined by power series and in general cannot be written in the standard form $A (X) = A \dot X \dot A^*$}.
			\item The orbit 
			$$C = C(\cn) := G(I) = \{X = A A^* \} \subset \Herm_m$$
			of the identity matrix $I$ is a homogeneous convex cone and the group $G$ acts on $C$ freely ( i.e. with trivial stabilizer $G_I = \{e\}).$
			\item Any homogeneous convex cone $C \subset V \simeq \bR^N$ is linearly isomorphic to the cone $C(\cn) \subset \Herm(\cn)$ associated to a $\Nil$-algebra $\cn$, defined up to an equivalence (see Definition \ref{Vinberg}).
			\item The adjoint cone $C^*= \{\xi \in W^* , \xi(x)>0, \, \forall x \in C \} \subset W^* = \Hom(W, \bR)$
			is isomorphic to the cone $C(\cn^{\tau})$, associated with the dual (anti-transposed) Nil-algebra $\cn^{\tau}$.
		\end{itemize}
		\et

		\subsection{Classification of self-adjoint cones}
		\bt (Jordan-Von Neumann-Wigner-K\"ocher-Vinberg)
		Let $\cn$ be a rank $n$ $\Nil$ algebra and $C(\cn)\subset \cH = \Herm(\cn) $ the associated homogeneous convex cone in the space $\Herm(\cn)$ of Hermitian matrices.
		The following conditions are equivalent:
		\begin{itemize}
			\item[i)] The algebra $\cn$ is indecomposable and isomorphic to the dual algebra $\cn^{\tau}$.
			\item[ii)] The algebra $\Herm_n(\cn)$ of Hermitian matrices is
a Euclidean indecomposable Jordan algebra, that is, the algebra $\Herm_n(\bK)$ of matrices over the division algebra $\bK = \bR,\bC, \bH$ and of $\bO$ for $n=3$, or, when $n=2$ the spin factor algebra $\Herm_2(V)$.
			\item[iii)] The cone $C(\cn)$ is indecomposable and self-adjoint.
 \item[iv)] The cone $C(\cn)$
 is indecomposable and is a symmetric space with respect to the canonical metric. 
		\end{itemize}
	 \et

		\section{A description of (quasi)Nil-algebras in therms of directed acyclic graphs}

	\subsection{Directed acyclic graphs, associated with a quasiNil-algebras}

 \begin{Def}
 We say that a {\bf Nil-graph} is a directed acyclic $\Gamma$ graph of diameter $d(\Gamma)=1$. 
 \end{Def} 
 
		We associate with an indecomposable quasiNil-algebra $\cN$, defined up to an equivalency, Nil-graph $\G(\cN)$, called the { adjacency graph} of $\cN$. QuasiNil-algebras $\cN$ with given graph $\G$ bijectively correspond to the equipment of the graph $\G$, see below.
\begin{Def}
 Let $\cN=\oplus_{i<j} \cN_{ij}$ be a quasiNil-algebra.
 \begin{itemize}
 \item [i)] The matrix $A = A(\cn)=||a_{ij}||$, where
$$
a_{ij}=\begin{cases}
 1 \ \ \ \text{if} \ \ \ \cN_{ij}\ne \{0\}; \\
 0 \ \ \ \text{if} \ \ \ \cN_{ij}=\{0\};
\end{cases}
$$
 is called the {\bf adjacency matrix} of $\cN$. 

 \item[ii)] The graph $\G(\cN) :=\Gamma(A(\cN))$ with the adjacency matrix $A(\cN)$ is called the {\bf adjacency graph} of $\cN$.

 \item [iii)] The length of a path in a graph $\G$ is the number of its edges. We denote by $ML(\G)$ the maximal length of paths in $\G$.

 \end{itemize}

\end{Def}

\bp \begin{itemize} \label{graphs} \ 
 \item [i)] Equivalent (quasi)Nil-algebras $\cN'$ and $\cN''$ have the same adjacency graphs.
 \item [ii)] Any Nil-graph $\G$ of diameter 1 is the adjacency graph of some Nil-algebra $\cN$.
 \item [iii)] The nilpotency index $NI(\cN)$ of a (quasi)Nil algebra $\cN$ is equal to the maximal length $ML(\G(\cN))$ of paths.
\end{itemize}
\ep
\begin{proof}
 i) This item is obvious. An admissible permutation $\sigma \in S_m$ corresponds to a renumbering of vertices. 
 
 ii) We may assume that the adjacency matrix $A = ||A_{ab}||$ of the graph $\G$ is nilpotent and upper triangular. Changing all non-zero entries $A_{ij} =1$ to a one of real associative division algebras $\bR,\, \bC,\, \bH$, we get a Nil-algebra $\cN$.
 
iii) Let $\G= (P, E)= ((1. \cdots.n) ,(e_{ij})_{t<j})$ be the graph of $\cN$ and 
$ e_{i_1 i_2}\ldots e_{i_{d}i_{d+1}}$ a maximal path in the graph $\G$. It corresponds to the maximal non-zero product 
$$\cN_{i_1 i_2}\cdot \cN_{i_2 i_3}\cdot \cdots \cdot \cN_{i_d i_{d+1}}$$
in the algebra $\mcal{N}$. Hence, $NI(\cN) = ML(\G(\cN)).$ 
 \end{proof}

\begin{Def}\label{defequipments}
Let $\Gamma$ be a Nil-graph directed acyclic graph of diameter $d(\Gamma)=1$. An {\bf equipment} of $\Gamma$ is a family of Euclidean spaces $\{\cN_{ab}\}$ defined for all edges $ab\in E$ and a family of isometric maps $\mu_{abc}:\cN_{ab}\times \cN_{bc} \to\cN_{ac}$ associated with all paths of length 2 $abc\in E$. We say that an equipment on a graph $\G$ 
 is {\bf Clifford}, (resp. {\bf associative} , resp., {\bf admissible}) if the corresponding
 quasiNil-algebras (see the construction below) are Clifford, (resp.,associative , resp. Nil-algebras).
 \end{Def}
 \begin{Def}
 Let $(\{\cN_{ab}\},\{\mu_{abc}\})$ be an equipment of a Nil-graph $\G$ of diameter 
 ${d(\Gamma)=1}$ and $i:\G \to \{1,\ldots, n\}$ a numbering such that $i(a)<i(b)$ for any $ab\in E$. Then the 
 numbering $i$ associates with the equipment $(\{\cN_{ab}\},\{\mu_{abc}\})$ of $\G$ the quasiNil-algebra $\oplus \cN_{i(a)i(b)}:=\cN_{ab}$ by $(\{\cN_{ab}\},\{\mu_{abc}\})$. The algebra $\oplus \cN_{i(a)i(b)}:=\cN_{ab}$ by $(\{\cN_{ab}\},\{\mu_{abc}\})$ is called {\bf an associated with the equipment $(\{\cN_{ab}\},\{\mu_{abc}\})$ of graph $\G$ quasiNil-algebra}.
\end{Def}

An associated quasiNil-algebra depends on the choice of a numbering, but its equivalence class does not depend.

 Thus, we have the following.
 \begin{Prop}
 There is a 1-1 correspondence between equivalence classes of indecomposable quasiNil-algebras and (connected) Nil-graphs with equipments.
 \end{Prop}
We propose the following strategy for the classification of rank $n$ Clifford Nil-algebras, that
 consists of three steps. 
 \begin{enumerate}
 \item [i)] Enumerate (connected) Nil-graphs with $n$ vertices of diameter 1.
 \item [ii)] Describe all Clifford equipments of these graphs $\G$, and associated quasiNil-algebras.
 \item[iii)] Describe when the constructed quasiNil-algebras are Nil-algebra i.e. when the algebra is associative and satisfies the Vinberg condition \ref{Vinberg condition} .
 \end{enumerate}
 Note that the first two steps reduce the problem of classification of all Clifford quasiNil-algebras to a purely combinatoric problem. 
 
 Recall that the nilpotency index $NI(\cN)$ of a quasiNil-algebra $ \cN$ is equal to the maximal length of the paths of the graph $\G(\cN)$ and a quasiNil-algebra $\cN$ with $NI(\cN)\le 2 $
 is a Nil-algebra (see \ref{rank4vinbercondition}).
 
 The Clifford Nil-algebras of rank 4 have nilpotency index $NI =1,2,3$. The Clifford Nil-algebras with $NI \leq 2$ are classified in Theorem \ref{index2} and Clifford quasiNil-algebras with maximal index $NI=3$ in Theorem \ref{cN-4}.
 
 The description of Clifford Nil-algebras with $NI =3$ is based on the new notion of Clifford extension (Definition \ref{DCEidentity}) and on $\bZ/2\bZ$-bigraded modules over a tensor product of two Clifford algebras $C_p\otimes C_q$. The classification of Clifford extensions is given in Theorem \ref{DCEclassification}. Bigraded modules over $C_p\otimes C_q$ are described in Section \ref{sectionbimodule}. 
 The classification of rank 4 Clifford Nil-algebras with nilpotency index $NI =3$ is given in Theorem \ref{index3}
 Below, we present the result of classification of rank 4 Clifford Nil-algebras 
 \subsection{Classification of Clifford quasiNil-algebras of rank 4}

 For brevity, we list Nil-algebras up to duality (i.e. reflection across the anti-diagonal). The corresponding graphs coincide up to inversion of all arrows. The equivalence class of a Nil-algebra does not depend on a numbering of vertices. For any graph, we choose a numbering, and it does not matter which numbering we choose.

		\subsubsection {Classification of rank 4 Clifford Nil-algebras $\cN$ with $ NI(\cN)\le 2 $}

 Recall that any quasiNil-algebra $\cN$ of rank 4 and nilpotency index $NI(\cN)\le 2$ is a Nil-algebra (Remark \ref{rank4vinbercondition}). That is, any equipment on a Nil-graph of diameter 1 defines a Nil-algebra.

If a Nil-algebra $\cN$ has the nilpotent index $NI(\cN) =1$, then all multiplications are trivial, and the equipment consists of associating a Euclidean vector space with each edge. 

 The quasiNil-algebras with the adjacency graph $ \G$ are Nil-algebras and correspond to equipments of edges $(ij)\in E$
 with Euclidean vector spaces $\cN_{ij}$. The following table lists such Nil-graphs and the associated quasiNil-algebras up to equivalence.

 \medskip	
		\begin{tabular}{|c|c|}
			\hline
			\begin{tikzpicture}
				\tikzset{vertex/.style = {circle,draw,fill=black,minimum size=0.5em, inner sep = 0pt}}
				\tikzset{edge/.style = {->,> = latex',thick}}
				\node[vertex][label=$1$] (a) at (0,0) {};
				\node[vertex][label=$2$] (b) at (2,0.8) {};
				\node[vertex][label=$3$] (c) at (2,0) {};
				\node[vertex][label=$4$] (d) at (2,-0.8) {};
				\draw[edge] (a) to (b);
				\draw[edge] (a) to (c);
				\draw[edge] (a) to (d);
			\end{tikzpicture} 
			& 
			{	 \begin{tikzpicture}
					\tikzset{vertex/.style = {minimum size=0.5em, inner sep = 0pt}}
					\node[vertex](b) at (0,0)
					{
						$
						\begin{pmatrix}
							0 & V_1 & V_2 & V_3 \\
							0 & 0 & 0 & 0 \\
							0 & 0 & 0 & 0 \\
							0 & 0 & 0 & 0
						\end{pmatrix}
						$};
			\end{tikzpicture} }
			\\ 
			
			\hline
			\begin{tikzpicture}
				\tikzset{vertex/.style = {circle,draw,fill=black,minimum size=0.5em, inner sep = 0pt}}
				\tikzset{edge/.style = {->,> = latex',thick}}
				\node[vertex][label=$1$] (a) at (0,1) {};
				\node[vertex][label=$2$] (b) at (0,0) {};
				\node[vertex][label=$3$] (c) at (2,1) {};
				\node[vertex][label=$4$] (d) at (2,0) {};
				\draw[edge] (a) to (c);
				\draw[edge] (a) to (d);
				\draw[edge] (b) to (c);
				\draw[edge] (b) to (d);
				\filldraw[black] (0,-0.4) circle (0pt);
			\end{tikzpicture} & 
			{	 \begin{tikzpicture}
					\tikzset{vertex/.style = {minimum size=0.5em, inner sep = 0pt}}
					\node[vertex](b) at (0,0)
					{
						$
						\begin{pmatrix}
							0 & 0 & V_1 & V_2 \\
							0 & 0 & V_3 & V_4 \\
							0 & 0 & 0 & 0 \\
							0 & 0 & 0 & 0
						\end{pmatrix}
						$};
			\end{tikzpicture} }
			\\
			\hline
			\begin{tikzpicture}
				\tikzset{vertex/.style = {circle,draw,fill=black,minimum size=0.5em, inner sep = 0pt}}
				\tikzset{edge/.style = {->,> = latex',thick}}
				\node[vertex][label=$1$] (a) at (0,1) {};
				\node[vertex][label=$2$] (b) at (0,0) {};
				\node[vertex][label=$3$] (c) at (2,1) {};
				\node[vertex][label=$4$] (d) at (2,0) {};
				\draw[edge] (a) to (c);
				\draw[edge] (a) to (d);
				\draw[edge] (b) to (d);
				\filldraw[black] (0,-0.4) circle (0pt);
			\end{tikzpicture} 
			& 	
			\begin{tikzpicture}
				\tikzset{vertex/.style = {minimum size=0.5em, inner sep = 0pt}}
				\node[vertex](b) at (0,0)
				{
					$
					\begin{pmatrix}
						0 & 0 & V_{1} & V_{2} \\
						0 & 0 & 0 & V_3 \\
						0 & 0 & 0 & 0 \\
						0 & 0 & 0 & 0
					\end{pmatrix}
					$};
			\end{tikzpicture} 
			\\
			\hline
		\end{tabular}
 
 \medskip

	 Similarly, if a Nil-algebra $\cN$ has the nilpotent index $NI(\cN) =2$, a quasiNil-algebra defined by an equipment of a connected Nil-graph $\G$ with $ML(\G) =2$ , is a Nil-algebra. All such Nil-graphs (up to a duality) and the restriction on the associated quasiNil-algebras are drawn below.

		\begin{tabular}{|c|c| }
			\hline
			\begin{tikzpicture}
				\tikzset{vertex/.style = {circle,draw,fill=black,minimum size=0.5em, inner sep = 0pt}}
				\tikzset{Vertex/.style = {minimum size=0.5em, inner sep = 0pt}}
				\tikzset{edge/.style = {->,> = latex',thick}}
				\node[vertex][label=$1$](a) at (0,0){};
				\node[vertex][label=$2$] (b) at (1.4,0) {};
				\node[vertex][label=$4$] (c) at (2.8,0) {};
				\node[vertex][label=$3$] (d) at (1.4,-1) {};
				\node[Vertex](f) at (0,-1.4) {};
				\draw[edge] (a) to (b);
				\draw[edge] (b) to (c);
				\draw[edge] (a) to (d);
				\draw[edge] (a) to[bend left] (c);
			\end{tikzpicture} & 
			{ \begin{tikzpicture}
					\tikzset{vertex/.style = {minimum size=0.5em, inner sep = 0pt}}
					\node[vertex](c) at (0,0.2) {};
					\node[vertex](b) at (0,-1){$
						\begin{pmatrix}
							0 & \cn_{1 2} & \cn_{1 3} & \cn_{1 4} \\
							0 & 0 & 0 &\cn_{2 4} \\
							0 & 0 & 0 & 0 \\
							0 & 0 & 0 & 0
						\end{pmatrix}$};
				\end{tikzpicture} 
				
			} \\
			\hline 
			\begin{tikzpicture}
				\tikzset{vertex/.style = {circle,draw,fill=black,minimum size=0.5em, inner sep = 0pt}}
				\tikzset{edge/.style = {->,> = latex',thick}}
				\node[vertex][label=$3$] (c) at (2.8,1.2) {};
				\node[vertex][label=$4$] (d) at (2.8,-1.2) {};
				\node[vertex][label=$1$] (a) at (0,0) {};
				\node[vertex][label=$2$] (b) at (1.4,0) {};
				\draw[edge] (a) to (b);
				\draw[edge] (a) to (d);
				\draw[edge] (a) to (c);
				\draw[edge] (b) to (c);
				\draw[edge] (b) to (d);
			\end{tikzpicture} 
			&
			{ \begin{tikzpicture}
					\tikzset{vertex/.style = {minimum size=0.5em, inner sep = 0pt}}
					\node[vertex](c) at (0,-1.2) {};
					\node[vertex](b) at (0,0){$
						\begin{pmatrix}
							0 & \cn_{1 2} & \cn_{1 3} & \cn_{1 4} \\
							0 & 0 & \cn_{2 3} & \cn_{2 4} \\
							0 & 0 & 0 & 0\\
							0 & 0 & 0 & 0
						\end{pmatrix} $};
				\end{tikzpicture} 
				
			} \\ 
			\hline
			\begin{tikzpicture}
				\tikzset{vertex/.style = {circle,draw,fill=black,minimum size=0.5em, inner sep = 0pt}}
				\tikzset{edge/.style = {->,> = latex',thick}}
				\node[vertex][label=$2$] (b) at (1.4,1.2) {};
				\node[vertex][label=$3$] (c) at (1.4,-1.2) {};
				\node[vertex][label=$1$] (a) at (0,0) {};
				\node[vertex][label=$4$] (d) at (2.8,0) {};
				\draw[edge] (a) to (b);
				\draw[edge] (a) to (d);
				\draw[edge] (a) to (c);
				\draw[edge] (c) to (d);
				\draw[edge] (b) to (d);
			\end{tikzpicture} &
			{ \begin{tikzpicture}
					\tikzset{vertex/.style = {minimum size=0.5em, inner sep = 0pt}}
					\node[vertex](c) at (0,-1.2) {};
					\node[vertex](b) at (0,0){$
						\begin{pmatrix}
							0 & \cn_{1 2} & \cn_{1 3} & \cn_{1 4} \\
							0 & 0 & 0 &\cn_{2 4} \\
							0 & 0 & 0 & \cn_{3 4} \\
							0 & 0 & 0 & 0
						\end{pmatrix} $};
				\end{tikzpicture} 
				
			} \\
			\hline
		\end{tabular}

To describe all equipments, we will use the following remark.

 \begin{rem}
 	For any $1\le i <j\le 4 $ denote $\nu_{ij}= \dim \cN_{ij}$.
	 Suppose $1\le i<j<k\le 4 $ and $0<\nu(i,j)\le \nu (j,k)$ (or $\nu(i,j)\ge \nu (j,k)>0$). Then, according to Theorem \ref{ABSh}, there exists a vector space $V$ and a $\Cl(V)$-module $S^0\oplus S^1$ such that the ordered triple $(\cn_{ij}, \cn_{jk},C_{ik})$ with the isometric map $\varphi_{ijk}:{\cn_{ij}\times \cn_{jk} \to \cn_{ik}}$ is isomorphic to the ordered pair $(V, S^0, S^1)$ with the {isometric Clifford map} ${V\times S^0 \to S^1}$ (or $(S^0,V,S^1)$ with the isometric Clifford map ${S^0\times V \to S^1}$). 	
 Thus, for any $1\le i<j<k\le 4$ that satisfies $\nu_{ij},\nu_{jk}\ne 0$ we have two possible cases.
 
 \end{rem}
		\begin{tabular}{|c|c|c|c| }
			\cline{1-3}
			\ 
			&
			{$\nu_{ij}\le\nu_{jk}$
				
			} & {$\nu_{ij}\ge\nu_{jk}$} \\
			\cline{2-3}
			\begin{tikzpicture}
				\tikzset{vertex/.style = {circle,draw,fill=black,minimum size=0.5em, inner sep = 0pt}}
				\tikzset{Vertex/.style = {minimum size=0.5em, inner sep = 0pt}}
				\tikzset{edge/.style = {->,> = latex',thick}}
				\node[vertex][label=$i$](a) at (0,0){};
				\node[vertex][label=$j$] (b) at (1.6,0) {};
				\node[vertex][label=$k$] (c) at (2.2,-1.4) {};
				\node[Vertex](f) at (0,-1.4) {};
				\draw[edge] (a) to (b);
				\draw[edge] (b) to (c);
				\draw[edge] (a) to (c);
			\end{tikzpicture} &
			{ \begin{tikzpicture}
					\tikzset{vertex/.style = {minimum size=0.5em, inner sep = 0pt}}
					\node[vertex](c) at (0,0.2) {};
					\node[vertex](b) at (0,-1){$
						\begin{pmatrix}
							0 & V & S^1 \\
							0 & 0 &S^0 \\
							0 & 0 & 0\\
						\end{pmatrix}$};
				\end{tikzpicture} 
				
			} & {	\begin{tikzpicture}
					\tikzset{vertex/.style = {minimum size=0.5em, inner sep = 0pt}}
					\node[vertex](c) at (0,0.2) {};
					\node[vertex](b) at (0,-1){$
						\begin{pmatrix}
							0 & S^0 & S^1 \\
							0 & 0 & V \\
							0 & 0 & 0\\
						\end{pmatrix}$};
				\end{tikzpicture} 
			} \\
			\cline{1-3}
		\end{tabular} \\ \\
		Using Remark 4.1, and considering all inequalities on $\nu_{ij},\nu_{jk}$ for all triangles $ijk$, we get all equipments given by the following notations: 
 $W, V, V_1, V_2$ are nonzero Euclidean spaces, $S^0\oplus S^1$ is a $\Cl(V)$-module, $T^0\oplus T^1$ is a $\Cl(S^1)$-module, $S^0_1\oplus S^1_1$ is a $\Cl(V_1)$-module, and $S^0_2\oplus S^1_2$ is a $\Cl(V_2)$-module. 
 \begin{center}
Table 1
\end{center}
		\begin{tabular}{|c|c|c|c| }
			\cline{1-3}
			\ 
			&
			{$\nu_{12}\le\nu_{13}$
				
			} & {$\nu_{12}\ge\nu_{13}$} \\
			\cline{2-3}
			\begin{tikzpicture}
				\tikzset{vertex/.style = {circle,draw,fill=black,minimum size=0.5em, inner sep = 0pt}}
				\tikzset{Vertex/.style = {minimum size=0.5em, inner sep = 0pt}}
				\tikzset{edge/.style = {->,> = latex',thick}}
				\node[vertex][label=$1$](a) at (0,0){};
				\node[vertex][label=$2$] (b) at (1.4,0) {};
				\node[vertex][label=$4$] (c) at (2.8,0) {};
				\node[vertex][label=$3$] (d) at (1.4,-1) {};
				\node[Vertex](f) at (0,-1.4) {};
				\draw[edge] (a) to (b);
				\draw[edge] (b) to (c);
				\draw[edge] (a) to (d);
				\draw[edge] (a) to[bend left] (c);
			\end{tikzpicture} &
			{ \begin{tikzpicture}
					\tikzset{vertex/.style = {minimum size=0.5em, inner sep = 0pt}}
					\node[vertex](c) at (0,0.2) {};
					\node[vertex](b) at (0,-1){$
						\begin{pmatrix}
							0 & V & W & S^1 \\
							0 & 0 & 0 & S^0 \\
							0 & 0 & 0 & 0 \\
							0 & 0 & 0 & 0
						\end{pmatrix}$};
				\end{tikzpicture} 
				
			} & {	\begin{tikzpicture}
					\tikzset{vertex/.style = {minimum size=0.5em, inner sep = 0pt}}
					\node[vertex](c) at (0,0.2) {};
					\node[vertex](b) at (0,-1){$
						\begin{pmatrix}
							0 & S^0 & W & S^1 \\
							0 & 0 & 0 & V \\
							0 & 0 & 0 & 0 \\
							0 & 0 & 0 & 0
						\end{pmatrix} $};
				\end{tikzpicture} 
			} \\
			\hline
			\ 
			&
			{$\nu_{23},\nu_{24}\le\nu_{12}$
				
			} & {$\nu_{23},\nu_{24}\ge\nu_{12}$}
			& {$\nu_{24}\le\nu_{12}\le \nu_{23}$} \\
			\cline{2-4}
			\begin{tikzpicture}
				\tikzset{vertex/.style = {circle,draw,fill=black,minimum size=0.5em, inner sep = 0pt}}
				\tikzset{edge/.style = {->,> = latex',thick}}
				\node[vertex][label=$3$] (c) at (2.8,1.2) {};
				\node[vertex][label=$4$] (d) at (2.8,-1.2) {};
				\node[vertex][label=$1$] (a) at (0,0) {};
				\node[vertex][label=$2$] (b) at (1.4,0) {};
				\draw[edge] (a) to (b);
				\draw[edge] (a) to (d);
				\draw[edge] (a) to (c);
				\draw[edge] (b) to (c);
				\draw[edge] (b) to (d);
			\end{tikzpicture} & 
			{ \begin{tikzpicture}
					\tikzset{vertex/.style = {minimum size=0.5em, inner sep = 0pt}}
					\node[vertex](c) at (0,-1.2) {};
					\node[vertex](b) at (0,0){$
						\begin{pmatrix}
							0 & S^0_1\simeq S^0_2 & S^1_0 & S^1_0 \\
							0 & 0 & V_1 & V_2 \\
							0 & 0 & 0 & 0\\
							0 & 0 & 0 & 0
						\end{pmatrix} $};
				\end{tikzpicture} 
				
			} & {	\begin{tikzpicture}
					\tikzset{vertex/.style = {minimum size=0.5em, inner sep = 0pt}}
					\node[vertex](c) at (0,-1.2) {};
					\node[vertex](b) at (0,0){$
						\begin{pmatrix}
							0 & V & S^1_1 & S^1_2 \\
							0 & 0 & S^0_1 & S^0_2 \\
							0 & 0 & 0 & 0 \\
							0 & 0 & 0 & 0
						\end{pmatrix} $};
				\end{tikzpicture} 
			} 
			& {	\begin{tikzpicture}
					\tikzset{vertex/.style = {minimum size=0.5em, inner sep = 0pt}}
					\node[vertex](c) at (0,-1.2) {};
					\node[vertex](b) at (0,0){$
						\begin{pmatrix}
							0 & S^0 & T^1 & S^1 \\
							0 & 0 & T^0 &V \\
							0 & 0 & 0 & 0 \\
							0 & 0 & 0 & 0
						\end{pmatrix} $};
				\end{tikzpicture} 
			}
			\\
			\hline
		\end{tabular}
		\\
		\begin{tabular}{|c|c|c|c|}
			\hline
			\ 
			&
			{$\nu_{12}\le \nu_{24} \ \text{and} \ \nu_{13}\ge\nu_{34}$
				
			} & {$\nu_{12}\le\nu_{24}\ \text{and} \ \nu_{13}\le\nu_{34}$} \\
			\cline{2-3}
			\begin{tikzpicture}
				\tikzset{vertex/.style = {circle,draw,fill=black,minimum size=0.5em, inner sep = 0pt}}
				\tikzset{edge/.style = {->,> = latex',thick}}
				\node[vertex][label=$2$] (b) at (1.4,1.2) {};
				\node[vertex][label=$3$] (c) at (1.4,-1.2) {};
				\node[vertex][label=$1$] (a) at (0,0) {};
				\node[vertex][label=$4$] (d) at (2.8,0) {};
				\draw[edge] (a) to (b);
				\draw[edge] (a) to (d);
				\draw[edge] (a) to (c);
				\draw[edge] (c) to (d);
				\draw[edge] (b) to (d);
			\end{tikzpicture} & 
			{ \begin{tikzpicture}
					\tikzset{vertex/.style = {minimum size=0.5em, inner sep = 0pt}}
					\node[vertex](c) at (0,-1.2) {};
					\node[vertex](b) at (0,0){$
						\begin{pmatrix}
							0 & V_1 & S_2^0 & S_1^1\simeq S_2^1 \\
							0 & 0 & 0 &S_1^0 \\
							0 & 0 & 0 & V_2 \\
							0 & 0 & 0 & 0
						\end{pmatrix} $};
				\end{tikzpicture} 
				
			} & {	\begin{tikzpicture}
					\tikzset{vertex/.style = {minimum size=0.5em, inner sep = 0pt}}
					\node[vertex](c) at (0,-1.2) {};
					\node[vertex](b) at (0,0){$
						\begin{pmatrix}
							0 & V_1 & V_2 & S_1^1\simeq S_2^1 \\
							0 & 0 & 0 &S_1^0 \\
							0 & 0 & 0 & S_2^0 \\
							0 & 0 & 0 & 0
						\end{pmatrix} $};
				\end{tikzpicture} 
			} 
			\\
			\cline{2-3}
			\ 
			&
			{$\nu_{12}\ge \nu_{24}\ \text{and} \ \nu_{13}\ge\nu_{34}$
				
			} & {$\nu_{12}\ge\nu_{24}\ \text{and} \ \nu_{13}\le\nu_{34}$} \\
			\cline{2-3}
			\ & 
			{ \begin{tikzpicture}
					\tikzset{vertex/.style = {minimum size=0.5em, inner sep = 0pt}}
					\node[vertex](c) at (0,-1.2) {};
					\node[vertex](b) at (0,0){$
						\begin{pmatrix}
							0 & S_1^0 & S_2^0 & S_1^1\simeq S_2^1 \\
							0 & 0 & 0 & V_1 \\
							0 & 0 & 0 & V_2 \\
							0 & 0 & 0 & 0
						\end{pmatrix} $};
				\end{tikzpicture} 
				
			} & {	\begin{tikzpicture}
					\tikzset{vertex/.style = {minimum size=0.5em, inner sep = 0pt}}
					\node[vertex](c) at (0,-1.2) {};
					\node[vertex](b) at (0,0){$
						\begin{pmatrix}
							0 & S_1^0 & V_2 & S_1^1\simeq S_2^1 \\
							0 & 0 & 0 & V_1 \\
							0 & 0 & 0 & S_2^1 \\
							0 & 0 & 0 & 0
						\end{pmatrix} $};
				\end{tikzpicture} 
			} 
			\\ 
			\hline
		\end{tabular} \\ \\
 These prove Theorem \ref{index2}
		
		
 
	\subsubsection {Classification of rank 4 Clifford quasiNil-algebras $\cN$ with $ NI(\cN)= 3 $}

 In addition to the standard notation $S = S^0 +S^1$ for a $\mathbb{Z}/2\mathbb{Z}$ (with an admissible Euclidean metric), we also fix some notation for different Clifford $CL(S^{\epsilon}),\, \epsilon=0,1)$ graded $CL(V)$ modules. We fix the following notations.
 \begin{center}
     Table 2
 \end{center}
 $$
 \begin{tabular}{|c|c|}
 \cline{1-2}
 an algebra $A$& an $A$-module \\
 \cline{1-2}
 $\Cl(V)$ & $S^0\oplus S^1$ \\
 \cline{1-2}
 $\Cl(S^0)$ & $T_1\oplus T_2$ \\
 \cline{1-2}
 $\Cl(S^1)$
 & $T_1\oplus T_3$ \\
 \cline{1-2}
 $\Cl(V)$
 & $T_2\oplus T_3$ \\
 \cline{1-2}
 $\Cl(V^{10})$
 & $S^{00}\oplus S^{10}$, \ $S^{10}\oplus S^{11}$ \\
 \cline{1-2}
 $\Cl(V^{01})$
 & $S^{01}\oplus S^{01}$, \ $S^{10}\oplus S^{11}$ \\
 \cline{1-2}
 \end{tabular} 
 $$

 There exists a unique Nil-graph on 4 vertices with the length of a maximal path 3 and diameter 1. Using Remark 4.1, and considering all inequalities on $\nu_{ij}=\cN_{ij},\nu_{jk}=\cN_{jk}$ for all triangles $ijk$, we get all equipments of this graph given by the following notations: $V,V^{10},V^{01}$ are Euclidean spaces and $T_1, T_2, T_3, S^{00},S^{10},S^{01},S^{11}$ as in the following table. 
\begin{center}
    Table 3
\end{center}
		\begin{tabular}{|c|c|c|}
			\hline
			\ 
			& 
			{$\nu_{12} \le \nu_{23}\le\nu_{34}$
				
			} & {$\nu_{34}\le\nu_{23}\le\nu_{12}$} \\
			\cline{2-3}
			\begin{tikzpicture}
				\tikzset{vertex/.style = {circle,draw,fill=black,minimum size=0.5em, inner sep = 0pt}}
				\tikzset{edge/.style = {->,> = latex',thick}}
				\node[vertex][label=$1$] (a) at (0,0) {};
				\node[vertex][label=$2$] (b) at (1.6,0) {};
				\node[vertex][label=$3$] (c) at (3.2,0) {};
				\node[vertex][label=$4$] (d) at (4.8,0) {};
				\draw[edge] (a) to (b);
				\draw[edge] (b) to (c);
				\draw[edge] (c) to (d);
				\draw[edge] (a) to[bend left] (c);
				\draw[edge] (b) to[bend right] (d);
				\draw[edge] (a) to[bend left] (d);
			\end{tikzpicture} 
			&
			{ \begin{tikzpicture}
					\tikzset{vertex/.style = {minimum size=0.5em, inner sep = 0pt}}
					\node[vertex](c) at (0,1) {};
					\node[vertex](b) at (0,0){$
						\cN_1^4:=\begin{pmatrix}
							0 & V & S^1 & T_3 \\
							0 & 0 & S^0 &T_2 \\
							0 & 0 & 0 & T_1 \\
							0 & 0 & 0 & 0
						\end{pmatrix} $};
				\end{tikzpicture} 
				
			} & {	\begin{tikzpicture}
					\tikzset{vertex/.style = {minimum size=0.5em, inner sep = 0pt}}
					\node[vertex](c) at (0,1) {};
					\node[vertex](b) at (0,0){$
						\left(\cN_1^4\right)^*:=\begin{pmatrix}
							0 & T_1 & T_2 & T_3 \\
							0 & 0 & S^0 & S^1 \\
							0 & 0 & 0 & V \\
							0 & 0 & 0 & 0
						\end{pmatrix} $};
				\end{tikzpicture} 
			} 
			\\
			\cline{2-3}
			\ 
			& 
			{$\nu_{23}\le \nu_{12}\le\nu_{34}$	
			} & {$\nu_{23}\le\nu_{34} \le\nu_{12}$} \\
			\cline{2-3}
			\ & 
			{ \begin{tikzpicture}
					\tikzset{vertex/.style = {minimum size=0.5em, inner sep = 0pt}}
					\node[vertex](c) at (0,1) {};
					\node[vertex](b) at (0,0){$
						\cN_2^4:=\begin{pmatrix}
							0 & S^0 & S^1 & T_3 \\
							0 & 0 & V & T_2 \\
							0 & 0 & 0 & T_1 \\
							0 & 0 & 0 & 0
						\end{pmatrix} $};
				\end{tikzpicture} 
			} & {	\begin{tikzpicture}
					\tikzset{vertex/.style = {minimum size=0.5em, inner sep = 0pt}}
					\node[vertex](c) at (0,1) {};
					\node[vertex](b) at (0,0){$
						\left(\cN_2^4\right)^*:=\begin{pmatrix}
							0 & T_1 & T_2 & T_3 \\
							0 & 0 & V & S^1 \\
							0 & 0 & 0 & S^0 \\
							0 & 0 & 0 & 0
						\end{pmatrix} $};
				\end{tikzpicture} 
			} 
			\\ 
			\cline{2-3}
			\ & {$\nu_{12},\nu_{34}\le\nu_{23}$ } \\
			\cline{2-2} 
			\ & {	\begin{tikzpicture}
					\tikzset{vertex/.style = {minimum size=0.5em, inner sep = 0pt}}
					\node[vertex](c) at (0,1) {};
					\node[vertex](b) at (0,0){$
						\cN_0^4:=\begin{pmatrix}
							0 & V^{10} & S^{10} & S^{11} \\
							0 & 0 & S^{00} & S^{01} \\
							0 & 0 & 0 & V^{01} \\
							0 & 0 & 0 & 0
						\end{pmatrix} $};
			\end{tikzpicture} }\\
			\cline{1-2} 
		\end{tabular}
		\\ \\
 
 Here, we assume that the algebras $\cN_{i}^4=\cN_{i}^4(V,S^0,S^1,T_1,T_2,T_3)$ depend on the collection $(V,S^0,S^1,T_1,T_2,T_3)$.
 
 \begin{theorem}\label{cN-4}
 Any Clifford quasiNil-algebra of rank 4 and nilpotency index 3 is isomorphic to $\cN_1^4,\cN_2^4$ or $\cN_0^4$, up to duality.
 \end{theorem}

 A quasiNil-algebra $\cN$ of rank 4 is a Nil-algebra if and only if $\cN$ is associative and satisfies Condition \eqref{Vinberg condition}. In Sections \ref{sectionCliffordExtension} we study quasiNil-algebras of the first two forms from Theorem \ref{cN-4} and find under what conditions these quasiNil-algebras are associative. In Section \ref{sectionbimodule}, we do the same for quasiNil-algebras of the last form from Theorem \ref{cN-4}. In Section \ref{sectionfinal}, we find which associative Clifford quasiNil-algebras are Nil-algebras.

 		\section{Clifford extensions} \label{sectionCliffordExtension}

 \subsection{Clifford quasiNil-algebras of nilpotency index 3 and Clifford extensions}

 In this section, we study algebras $\cN_1^4$ and $\cN_2^4$. Namely, we describe collections $$
 (V,S^0,S^1,T_1,T_2,T_3)
 $$
 satisfying the next definition.
		
	\begin{Def}
	 Let $S=S^0\oplus S^1$ be a $ \bZ/2\bZ$-graded module over $\Cl(V)$; $T_1,T_2,T_3$  Euclidean spaces of the same dimension, $T_2\oplus T_3$, $T_1\oplus T_2$, and $T_1\oplus T_3$  graded $\Z/2\Z$-modules over $\Cl(V),\Cl(S^0)$, and $\Cl(S^1)$, respectively (as in Table 2). Then we say that $(T_1,T_2,T_3)$ is a {\bf Clifford extension} over $(V,S^0,S^1)$.
	\end{Def}
				The associativity condition for the algebra $\cN_{1}^4$ reads as: 
 for any $v\in V, s^0\in S^0,t^0\in T^0$ the identity 
		$$
		(vs^0)t^0=v(s^0t^0).
		$$
		holds. Similarly, 	the associativity condition for the algebra $\cN_{2}^4$ 
 means that for any $v\in V, s^0\in S^0,t^0\in T^0$ the identity 
		$$
		(vs^1)t_2=s^1(vt_2).
		$$
 holds. 
		The described construction motivates the next definition. 
		
		\begin{Def}\label{DefDCE} 

 Let $T_1,T_2,T_3$ be a Clifford extension over $(V,S^0,S^1)$. Suppose that for any $v\in V, s^0\in S^0, s^1=v s^0 \in S^1$, the multiplication maps 
			$$
			\mu_v: T_2 \to T_3, \ \ \mu_{s^0}: T_1 \to T_2, \ \ \mu_{s^1}: T_1 \to T_3
			$$
			(that are isomorphisms) form a commutative diagram
			\begin{equation}\label{diagram}
				\xymatrix{
					&T_2\ar@<1ex>[dr]^{\mu_v}
					\ar@<0ex>[dl]\\
					T_1\ar@<1ex>[ur]^{\mu_{s^0}} \ar@<0ex>[rr]_{\mu_{s^1}}&&T_3\ar@<0ex>[ul] \ar@<-0.9ex>[ll]
				}.
			\end{equation}
			Then we say $(T_1,T_2,T_3)$ is an {\bf associative Clifford extension} over $(V,S^0,S^1)$. 
		\end{Def}

		\begin{Prop}\label{dCe}
			Let $T_1,T_2,T_3$ be a Clifford extension over $(V,S^0,S^1)$. Then the following conditions are equivalent:
			\begin{itemize}
				\item [i)] $(T_1,T_2,T_3)$ is an associative Clifford extension over $(V,S^0,S^1)$;
				\item[ii)] For any $v\in V, s^1\in S^1, t_1\in T_1$ we have $ (vs^1) t_1=v (s^1 t_1)$;
				\item[iii)] For any $v\in V, s^1\in S^1, t_2\in T_2$ we have $ (vs^1) t_2=s^1 (v t_2)$;
				\item[iv)] For any $v\in V, s^0\in S^0, t_2\in T_2$ we have $ vt_2=-(vs^0)(s^0t_2)$;
			\end{itemize} 
		\end{Prop}
		
		\begin{proof}

			For any $v\in V, s^0\in S^0, s^1=v s^0 \in S^1$, the maps 
			$$
			\mu_v: T_2 \to T_3, \ \ \mu_{s^0}: T_1 \to T_2, \ \ \mu_{s^1}: T_1 \to T_3
			$$
			are isomorphisms. 
			If one of the conditions
			$$
			1)\ \mu_{s^1}=\mu_v \mu_{s^0}, \ \ \ 2)\ \mu_{s^0}=\mu_{v}^{-1}\mu_{s^1}, \ \ \ 3)\ \mu_{v}=\mu_{s^1}\mu_{s^0}^{-1}
			$$
			is satisfied then the diagram \ref{diagram}
			is commutative. Moreover, if this diagram is commutative then the conditions 1), 2), 3) are satisfied. Hence, if one of these items is satisfied, then the other two are satisfied, too. 
			
			We have
			\begin{align*}
				\mu_{s^1}=\mu_{v}\mu_{s^0} \ &\Leftrightarrow \ \forall t_1\in T^1 \ (vs^0) t_1=v (s^0 t_1); \\ 
				\mu_{s^0}=\mu_{v}^{-1}\mu_{s^1} \ &\Leftrightarrow \ \forall t_1\in T^1 \ (vs^1) t_1=v (s^1 t_1); \\ 
				\mu_{v}=\mu_{s^1}\mu_{s^0}^{-1} \ &\Leftrightarrow \ \forall t_2\in T^2 \ vt_2=-(vs^0)(s^0t_2); 
			\end{align*}
			Therefore, if any of the condition i) -- iv) is satisfied then the diagram \eqref{diagram} is commutative and all conditions i) -- iv) are satisfied.
		\end{proof}
		 As a corollary, we get the following.
 \bt \label{associtiveSQNA}
 There exists a 1-1-1 correspondence between Clifford extensions $(T_1,T_2,T_3)$ over $(V,S^0,S^1)$, quasiNil-algebras $\cN_{1}^4(V,S^0,S^1,T_1,T_2,T_3)$, and quasiNil-algebras $\cN_{2}^4(V,S^0,S^1,T_1,T_2,T_3)$
		 Moreover, Clifford extension is associative if and only if the corresponding quasiNil-algebras are associative.
		\et

 \subsection{Classification of associative Clifford extensions}

 		\subsubsection{Restrictions on dimensions}
		
		\begin{Prop}\label{p45}
			Let $(T_1,T_2,T_3)$ be an associative Clifford extension over $(V,S^0,S^1)$. Then $\dim V=1$ or $\dim V\le \dim S^0 \le 4$.
		\end{Prop}
		
		The idea of the proof is the following. If $\dim V>1$ and $\dim S^0>4$ then we construct four orthogonal vectors $e_1, e_2,e_3,e_4 \in S^0$ with the same action of the elements $e_1e_2, e_3 e_4\in \Cl(S^0)$ in the module $S^0\oplus S^1$. Then we prove that $e_1e_2 \in Cl(S^0)$ and $e_3 e_4 \in Cl(S^0)$ must be equal, which is impossible. 
		
		\begin{Lem}\label{DCEidentity}
			Let $(T_1, T_2, T_3)$ be an associative Clifford extension over $(V,S^0, S^1)$. Then for any ${v_1,v_2\in V,} {s^0\in S^0}$ the action of the element
			$$
			\left(v_1 v_2 s^0\right)s^0 \in \Cl^0(S^0)
			$$
			on $T_2$ does not depend on $s^0$.
		\end{Lem}
		\begin{proof} For any $t_2\in T_2$ we have
			$$
			((v_1 v_2 s^0)s^0)t_2= (v_1 (v_2 s^0))(s^0 t_2).
			$$
			Applying item ii) from Proposition \ref{dCe} to $v_1\in V, v_2 s^0\in S^1, s^0 t_2 \in T^1$, we get that 
			$$
			(v_1 (v_2 s^0))(s^0 t_2)=v_1((v_2 s^0)(s^0 t_2)).
			$$
			According to item iv) from Proposition \ref{dCe}, we find that for any $t_2\in T_2$ the element
			$$
			((v_1 v_2 s^0)s^0)t_2=v_1((v_2 s^0)(s^0 t_2))=-v_1 v_2 t_2 \in T_1
			$$
			does not depend on $s_0$. 
		\end{proof}
		
		\begin{Lem}\label{e1e2e3e4}
			Let $( T_1, T_2, T_3)$ be an associative Clifford extension over $(V,S^0,S^1)$ and $\dim V>1$, $\dim S^0\ge 4$. Then there exist orthogonal vectors $e_1,e_2,e_3,e_4\in T^0$ with the same actions of $$e_1 e_2, e_3 e_4\in \Cl^0(S^0)$$ on $T_2$.
		\end{Lem}

		\begin{proof}
			Choose any $s^0_1\in S^1, v_1,v_2\in V$, such that $v_1 \perp v_2$ and $s^0_1\in S^0 $. 
			Then $v_1v_2=-v_2v_1$ in $Cl(V)$ and
			$$
			g_S(v_1v_2s_1^0,s^0)=-g_S(v_2s_1^0,s_1^0)= 0$$
 since $v_2s_1^0 \in S^1 \perp s_1^0$.
			That is, $v_1 v_2 s_1 \perp s_2$. Since $\dim S^0>2$ there is $s_2^0\in S^0$ such that $s_2^0\perp s_1^0, s_2^0\perp v_1 v_2 s_0$. According to Lemma \ref{DCEidentity}, the action of 
			$$
			(v_1v_2)s_1^0 \ \ \ \text{and} \ \ \ (v_1 v_2) s_2^0
			$$ 
			on $T^2$ are the same. We have constructed four orthogonal vectors 
			$$
			e_1=v_1 v_2 s_1^0, \ \ \ \ e_2=s_1^0 , \ \ \ \ e_3=v_1 v_2 s_2^0, \ \ \ \ e_4=s_2^0 
			$$
			with the same multiplication of $e_1 e_2, e_3 e_4\in \Cl^0(S^0)$ on $T_2$.
		\end{proof}

		\begin{proof}[Proof of Proposition \ref{p45}]
			Suppose $\dim V>1$ and $\dim S^0>4$. According to Lemma \ref{e1e2e3e4}, there exist orthogonal vectors $e_1,e_2,e_3,e_4\in T^0$ with the same actions of $e_1 e_2, e_3 e_4\in \Cl^0(S^0)$ on $T_2$. Since $\dim S^0>4$, there exists a vector $e_0$ that is orthogonal to $e_1,e_2,e_3,e_4$. Any element $t_1\in T_1$ admits a form $t_1=e_0 t_2$. We have 
			$$
			e_1e_2t_1=e_1 e_2 e_0 t_2=e_0 e_1 e_2 t_2= e_0 e_3 e_4 t_2=e_3 e_4 e_0t_2= e_3 e_4 t_1.
			$$
			Hence, $e_1 e_2$ and $e_3 e_4$ induce the same multiplications on $T_1\oplus T_2$. Therefore,
			$$
			e_1 e_2=e_3 e_4
			$$ 
			but this is impossible for orthogonal vectors. 
			
			Therefore, $\dim V=1$ or $\dim S^0\le 4$. If $\dim S^0\le 4$, then 
			$\dim V \in \{1,2,3,4\}$.

		\end{proof}
		
		\subsubsection{$\dim V=1$}
		\begin{Example}
			Let $\tilde V$ be a Euclidean space and $\tilde S^0\oplus \tilde S^1$ a $\Cl(\tilde V)$-module. Consider the matrix algebra
			$$
			\cn=
			\begin{pmatrix}
				0 &\bR& \tilde V & \tilde S^1 \\
				0 & 0 & \tilde V & \tilde S^1 \\
				0 & 0 & 0 & \tilde S^0 \\
				0 & 0 & 0 & 0
			\end{pmatrix}
			$$
		with standard multiplications $\bR\times \tilde V \to \tilde V, \bR\times \tilde S^1\to \tilde S^1, V\times \tilde S^0 \to \tilde S^1$. Since for any ${a\in \bR,}{\tilde v\in \tilde V,}{\tilde s^0\in \tilde S^0}$ we have $a(\tilde v\tilde s^0)=(a\tilde v)\tilde s^0$, $\cn$ is an associative Clifford quasiNil-algebra and the collection $( \tilde S^0, \tilde S^1,\tilde S^1)$ is an associative Clifford extension over $(\bR, \tilde V,\tilde V)$.
		\end{Example}
		
 \begin{Prop}
			Let $(V,S^0,S^1,T_1,T_2,T_3)$ be an associative Clifford extension and $\dim V=1$. Then $(V,S^0,S^1,T_1,T_2,T_3)$ is isomorphic to $(\bR, \tilde V,\tilde V, \tilde S^0, \tilde S^1,\tilde S^1)$.
		\end{Prop}
		\begin{proof}
			The Clifford extension $(V,S^0,S^1,T_1,T_2,T_3)$ contains all multiplications from the matrix algebra 
			$$
			\begin{pmatrix}
				0 &V& S^1 & T_3 \\
				0 & 0 & S^0 & T_2 \\
				0 & 0 & 0 & T_1 \\
				0 & 0 & 0 & 0
			\end{pmatrix}
			$$
			(we think it is better to give this matrix here for easy reading).
			Identify $V\simeq \bR$ and denote $S^0\simeq \tilde V$. The space $T_1\oplus T_2$ is a $\Cl(\tilde V)$-module. Denote $T_1\simeq \tilde S^0$ and $T_2\simeq \tilde S^1$. Since we have an isometric bilinear map $\bR \times S^0 \to S^1$, the multiplication on $1\in \bR$ is an isometry $S^0 \to S^1$. Thus, we can identify $S^0\simeq S^1\simeq \tilde V$. Analogically, the multiplication on $1\in \bR$ is an isometry $T_2 \to T_3$. Thus, we can identify $T_2\simeq T_3\simeq \tilde V$. Since for any $v\in V, S^0 \in S^0, t_1\in T_1$ we have $v(s^0t_1)=(vs^0)t_1$ and the multiplications $V\times S^0 \to S^1, S^0\times T_2 \to T_3$ are identified with the multiplications $\bR \times \tilde V \to \tilde V, \bR \times \tilde S^1 \to \tilde S^1$, respectively, the multiplication $S^1\times T_1 \to T_3$ coincides with the standard multiplication $\tilde V \times \tilde S^0 \to \tilde S^1$ after identification $S^1\simeq \tilde V, T_1\simeq \tilde S^0,T_3\simeq \tilde S^1$.
		\end{proof}
		
		\begin{Cor}
			There is a 1-1 correspondence between associative Clifford extensions $(T_1,T_2,T_3)$ over $(V,S^0,S^1)$ with $\dim V=1$ and Clifford modules $(\tilde V,\tilde S^0\oplus \tilde S^1)$.
		\end{Cor}

 \subsubsection{$\dim V\ge 2$}

			\begin{Example}\label{dimV>1}
			 Let $\mathbb{A} = \bC$ or $\bH$, $V \subset \bA$ be a subspace. Then
 $S = \bA \oplus \bA$ is a $\bZ_2$ graded $Cl(V)$-module and
 $(T_1,T_2,T_3): = (\bA^n,\bA^n,\bA^n)$ is an associative Clifford extension. 
 It defines the following $r=4$ Nil-algebras $\cN$ with $ NI(\cN) =3$: 
			$$
			\begin{pmatrix}
				0 & V & \bA & \bA^n \\
				0 & 0 & \bA & \bA^n \\
				0 & 0 & 0 & \bA^n \\
				0 & 0 & 0 & 0
			\end{pmatrix}, \ \ \ 
			\begin{pmatrix}
				0 & \bA & \bA & \bA^n \\
				0 & 0 & V & \bA^n \\
				0 & 0 & 0 & \bA^n \\
				0 & 0 & 0 & 0
			\end{pmatrix}, 
			\ \ \ 
			\begin{pmatrix}
				0 & \bA^n & \bA^n & \bA^n \\
				0 & 0 & V & \bA \\
				0 & 0 & 0 & \bA \\
				0 & 0 & 0 & 0
			\end{pmatrix},
			\ \ \ 
			\begin{pmatrix}
				0 & \bA^n & \bA^n & \bA^n \\
				0 & 0 & \bA & \bA \\
				0 & 0 & 0 & V\\
				0 & 0 & 0 & 0
			\end{pmatrix}.
			$$		
			
		\end{Example}
 The following lemma is well known.
		\begin{Lem} \label{lemmamu}
		 Let $V$ be the oriented Euclidean vector space of dimension 2 (respectively, 4). Then
 the structures of division algebra $\mathbb{K} =\bC, \bH$ with unit vector $e_0$ on $V$ bijectively corresponds to the oriented orthonormal frames $(e_0,i)$ or, respectively, $(e_0, i,j,k)$. Any two such structures are conjugated by an automorphism of $\mathbb{K} = \bZ_2$ or $ SO_3$.
		\end{Lem} 

		\begin{Prop}\label{WAAAAA}
			Let $(T_1,T_2,T_3)$ be a Clifford extension over $(V,S^0,S^1)$ and 
			$$\dim S^0=\dim S^1 =\dim T_1= \dim T_2=\dim T_3=k\in\{2,4\}.$$ 
 Then it is isomorphic to the Clifford extension of the form 
 $(V, \bA, \bA, \bA, \bA ,\bA)$, where $\bA\ =\bC$ or $\bH$ and $V\subset \bA$ is a vector subspace with the standard multiplication.
		\end{Prop}
		
		\begin{proof}
			Consider the Nil-algebra
			$$
			\begin{pmatrix}
				0 & V_{12} & V_{13} & V_{14} \\
				0 & 0 & V_{23} &V_{24} \\
				0 & 0 & 0 & V_{34} \\
				0 & 0 & 0 & 0
			\end{pmatrix}
			=
			\begin{pmatrix}
				0 & V & S^1 & T_3 \\
				0 & 0 & S^0 &T_2 \\
				0 & 0 & 0 & T_1 \\
				0 & 0 & 0 & 0
			\end{pmatrix}
			$$ 
			We have $$
			\dim V_{13} = \dim V_{14} = \dim V_{23} =\dim V_{24} = \dim V_{34} =k\in \{2,4\}.
			$$
			Choose
			$$
			e_{12}\in V_{12}, \ \ \ e_{23}\in V_{23}, \ \ \ e_{34}\in V_{34},
			$$
			$$
			e_{13}=e_{12}e_{23}\in V_{13},\ \ \ e_{24}=e_{23} e_{34}\in V_{24}, \ \ \ \ e_{14}=e_{12}e_{23}e_{34}\in V_{14}.
			$$
			Then we have the commutative diagram. 
			\begin{equation}
				\xymatrix{
					V_{12} \ar@{^{(}->}@<0ex>[r]^{R_{e_{23}}}& V_{13} \ar@<0.4ex>[d] \ar@<0.4ex>[r]^{R_{e_{34}}}& V_{14} \ar@<0.4ex>[l] \ar@<0.4ex>[d]\\
					& V_{23} \ar@<0.4ex>[u]^{L_{e_{13}}} \ar@<0.4ex>[r]^{R_{e_{34}}} & V_{24} \ar@<0.4ex>[l] \ar@<0.4ex>[u]^{L_{e_{13}}} \ar@<0.4ex>[d]\\ 
					& & V_{34} \ar@<0.4ex>[u]^{L_{e_{23}}}
				},
			\end{equation}
			where $R_{e_{23}}, R_{e_{34}}, L_{e_{13}},L_{e_{23}}$ are right and left multiplications on $e_{23}, e_{34}, e_{13},$. Moreover, the maps from the diagram send elements $e_{12},e_{13},e_{13},e_{23},e_{24},e_{34}$ to each other. Thus, we can identify
			$$
			(V_{13},e_{13})\simeq (V_{14},e_{14})\simeq (V_{23},e_{23})\simeq (V_{24},e_{24})\simeq (V_{34},e_{34})\simeq (V,e)
			$$
			and $V_{12}$ with a subset of $V$. According to Lemma \ref{lemmaXYZ}, triplets $(V_{12},V_{13},V_{23})$, $(V_{12},V_{14},V_{24})$, $(V_{13},V_{14},V_{34})$, $(V_{23},V_{24},V_{34})$ sets four quaternionic structures of a division algebra $\bA$ on $V$ with the same identity element. Let $\mu_{123},\mu_{124},\mu_{134},\mu_{234}$ be the corresponding multiplications in $V$. According to Lemma \ref{lemmamu}, for any $i,j,k,l$ we have
			$$
			\forall x,y\in V \ : \ \mu_{i,j}(x,y)=\mu_{k,l} (x,y) \ \ \ \text{or} \ \ \ 
			\forall x,y\in V \ : \ \mu_{i,j}(x,y)=\mu_{k,l} (y,x).
			$$
			The condition
			$$
			\forall v_{12}\in V_{12}, v_{23}\in V_{23}, v_{34}\in V_{34} \ : \- (v_{12}v_{23})v_{34}= v_{12}(v_{23}v_{34})
			$$
			is equivalent to 
			\begin{equation}\label{e42}
				\forall x,y,z \in V \mu_{134}(\mu_{123}(x,y),z) =\mu_{124}(x,\mu_{234}(y,z)).
			\end{equation}
			Combining this with Lemma \ref{lemmamu} we get that, \eqref{e42} can be satisfied if and only if all multiplications $\mu_{ij}$ coincide.
		\end{proof}
		
		\begin{Prop}\label{WAAAAAn}
			Let $(T_1,T_2,T_3)$ be an associative Clifford extension over $(V,S^0,S^1)$ and $\dim S^0=\dim S^1 =\dim T_1= \dim T_2=\dim T_3=k\in\{2,4\}$ then $(V,S^0,S^1,T_1,T_2,T_3)$ is isomorphic to $(W, \bA, \bA, \bA^n, \bA^n ,\bA^n)$, where $n\in \N, \bA\in \{\bC,\bH\}$, $W\subset \bA$ is a vector subspace, and all multiplication coincides with the standard multiplication on elements from $\bA$. 
		\end{Prop}
		
		\begin{proof}
			Take $t_1\in T^1$. Remind that $T_1\oplus T_2, T_1\oplus T_3, T_2\oplus T_3$ are $\Cl(S^0),\Cl(S^1),\Cl(V)$ modules, respectively. Then 
			$$
			\left(\Cl(S^0)\right)t_1=T_1'\oplus T_2'\subset T_1\oplus T_2 \ \ \ \text{and} \ \ \ \left(\Cl(S^1)\right)t_1=T_1'\oplus T_3'\subset T_1\oplus T_3
			$$ 
			are indecomposable submodules. For any $t_2'\in T_2,v\in V$ we have 
			$$
			vt_2'=-\frac{1}{g(s^0,s^0)}v(s^0 (s^0 t_2'))=-\frac{1}{g(s^0,s^0)}(vs^0) (s^0 t_2')\in T_3'
			$$
			because $s_0t_2'\in T_1', vs_0\in S^1$.
			
			Analogically, for any $t_3'\in T_3'$ we have 
			$$
			vt_3'=-\frac{1}{g(s^1,s^1)}v(s^1 (s^1 t_3'))=-\frac{1}{g(s_1,s_1)}(vs^1) (s^1 t_3')\in T_2',
			$$
			where we use item ii) of Proposition \ref{dCe} for $v\in V, s^1\in S^1, s^1t_3'\in T_1$. 
			
			We have checked that $T_2'\oplus T_3'\subset T_2\oplus T_3'$ is a submodule. 
			
			Consider the orthogonal decomposition 
			$$
			T_1=T_1'\oplus T_1'', \ \ \ \ \ T_2=T_2'\oplus T_2'', \ \ \ \ \ T_1=T_3'\oplus T_3''.
			$$
			Since $T_2\oplus T_3$ is an admissible $\Cl(V)$-module, for any $v\in V, t_2''\in T'', t_3'\in T_3'$ we have 
			$$
			g(vt_2'',t_3')=-g(t_2'',vt_3')=0. 
			$$ 
			Hence, 
			$$
			vt_2''\in T_3''.
			$$ 
			Analogically, we get that for any $v\in V, t_3''\in T_3''$ we have 
			$$
			vt_3''\in T_2''.
			$$ 
			
			Therefore, $T_2''\oplus T_3''$ is a $\Cl(V)$-module. 
 By the same way, we can prove that $T_1''\oplus T_2''$ is a $\Cl(S^0)$-module and $T_1''\oplus T_3''$ is a $\Cl(S^1)$-module. Thus, $(V,S^0, S^1, T_1',T_2',T_3')$ and $(V,S^0, S^1, T_1'',T_2'',T_3'')$ are Clifford extensions. Using Proposition \ref{WAAAAA} and mathematical induction we get that $(V,S^0, S^1, T_1,T_2,T_3)$ is isomorphic to $(W,\bA,\bA,\bA^n,\bA^n,\bA^n)$.
		\end{proof}
		
	We have proved the following.

	\begin{theorem}\label{DCEclassification}
		Let $(T_1,T_2,T_3)$ be an associative Clifford extension over $(V,S^0,S^1)$. 
		\begin{itemize}
			\item[a)] If $\dim V=1$ then $(V,S^0,S^1,T_1,T_2,T_3)$ is isomorphic to the collection $(\bR, \tilde V, \tilde V, \tilde S^0, \tilde S^1,\tilde S^1)$
			unique defined by a Euclidean vector space $S^0$ and a $\Cl(S^0)$-module $T_1\oplus T_3$.
			\item[b)] If $\dim V>1,\dim S^0=2$ then $(V,S^0,S^1,T_1,T_2,T_3)$ is equivalent to $(\bC,\bC,\bC,\bC^n,\bC^n,\bC^n)$.
			
			\item[c)] If $\dim V>1,\dim S^0>2$ then $(V,S^0,S^1,T_1,T_2,T_3)$ is equivalent to $(W,\bH,\bH,\bH^n,\bH^n,\bH^n)$, where $W\subset \bH$ is a vector subspace as in Example \ref{dimV>1}.
		\end{itemize}
	\end{theorem}

		\section{Bigraded \texorpdfstring{$C_p\otimes C_q$}{}-modules}\label{sectionbimodule}
 \subsection{Nil-algebras constructed from bigraded \texorpdfstring{$C_p\otimes C_q$-modules}{}}
		In this section, we describe when a quasiNil-algebra of the form 
 \begin{equation*}\label{bimod}
 \cN_0^4=\begin{pmatrix}
			0 & V^{10} & S^{10} & S^{11} \\
			0 & 0 & S^{00} & S^{01} \\
			0 & 0 & 0 & V^{01} \\
			0 & 0 & 0 & 0
		\end{pmatrix}
 \end{equation*}
 is associative. 
 We recall the notations in $\cn_0^4$: $S^{00}\oplus S^{10}, S^{01}\oplus S^{11}$ are $\Cl(V^{10})$-modules, $S^{00}\oplus S^{01}$, $S^{10}\oplus S^{11}$ are $\Cl(V^{01})$-modules (see Table 2). 
 
 The associativity condition reads as: for any $v\in V^{10}, v\in V^{01}, s^{00}\in S^{00}$ 
		$$
		v^{10}(v^{01} s^{00})= v^{01}(v^{10}s^{00}),
		$$ 
		is satisfied. Therefore, for any $v \in V^{10}, v'\in V^{01}$ we have the commutative diagram
 $$
 	\xymatrix{
					S^{10} \ar@<0.4ex>[d] \ar@<0.4ex>[r]^{\mu_{v'}}& S^{11} \ar@<0.4ex>[l] \ar@<0.4ex>[d]\\
					S^{00} \ar@<0.4ex>[u]^{\mu_{v}} \ar@<0.4ex>[r]^{\mu_{v'}} & S^{01} \ar@<0.4ex>[l] \ar@<0.4ex>[u]^{\mu_{v}} 
				},
 $$
 where $\mu_v$ and $\mu_v'$ are multiplications on $v$ and $v'$. Hence, for any $k,l\in\{0,1\},v^{10}\in V^{10}, v^{01}\in V^{01}, s^{kl}\in S^{kl}$ we have 
 $$
 v^{10}(v^{01} s^{kl}):= v^{01}(v^{10}s^{kl})
 $$
 Therefore, the action of $\Cl (V^{10})$ and $\Cl(V^{01})$ on $S:=S^{00}\oplus S^{10}\oplus S^{01}\oplus S^{11}$ extends to the action of $\Cl(V^{10})\otimes \Cl(V^{01})$ on $S$.


		A tensor product $\Cl(V^{10}) \otimes \Cl(V^{01})$ is the algebra with generators $u_1,\ldots, u_p, w_1,\ldots w_q$ and relations 
		$$
		u_i^2=-1, \ \ u_i^2=-1, \ u_iu_j=-u_ju_i, \ w_iw_j=-w_jw_i, \ u_iw_j=w_ju_i,
		$$ 
		for any $i\ne j$. 
		
		Endow the algebra $\Cl(V^{10}) \otimes \Cl(V^{01})$ with a bigradation 
		by the rule 
		$$
		u_1, \ldots, u_p \in \left(\Cl(V^{10}) \otimes \Cl(V^{01})\right)^{10}, \ \ \ \ w_1, \ldots, w_q \in \left(\Cl(V^{10}) \otimes \Cl(V^{01})\right)^{01}. 
		$$ 
		Then $S=S^{00}\oplus S^{10}\oplus S^{01}\oplus S^{11}$ is a $ \bZ/2\bZ$-bigraded $\Cl(V^{10}) \otimes \Cl(V^{01})$-module.

 Conversely, we can construct a Clifford Nil-algebra of the form $\cN_0^4$ by any $\bZ/2\bZ$-bigraded $\Cl(V^{10}) \otimes \Cl(V^{01})$-module $S$.

		\begin{Prop}\label{bimodulealgebra}
			There is a correspondence between associative Clifford quasiNil-algebras of the form $\cN_0^4$
			and $\bZ/2\bZ$-bigraded $\Cl(V^{10}) \otimes \Cl(V^{01})$-modules.
		\end{Prop}
 In section 7, we check that any associative Clifford quasiNil-algebras of the form $\cN_0^4$ is a Nil-algebra.

 \begin{Example}
 For any $p\in \N$ denote $C_p=\Cl(\bR^p)$. A tensor product of a graded $C_p$-module and a $C_q$-module is a bigraded $C_p\otimes C_q$-module. 
 \end{Example}

 \subsection{Describing of algebras \texorpdfstring{$C_p\otimes C_q$}{} in matrix forms}
 
		\begin{Prop}\label{bimodule}
			Let $\R: S=S^{00}\oplus S^{10}\oplus S^{01}\oplus S^{11} \to S^{00}$ be a functor which assigns to a bigraded $C_p\otimes C_q$-module $S$ the $(C_p\otimes C_q)^{00}$-modules $S^{00}$. Then $\R$ induces a 1-1 correspondence between bigraded $C_p\otimes C_q$ modules and (not graded) $(C_p\otimes C_q)^{00}$-modules.
		\end{Prop}
		\begin{proof}
			For $S^{00}$ is a $(C_p\otimes C_q)^{00}$-module set
			$$
			{\cal Q}(S^{00})=S^{00}\otimes_{(C_{p}\otimes C_q)^{00}} (C_{p}\otimes C_q).
			$$ 
			Then $\cal Q$ is a functor between $(C_p\otimes C_q)^{00}$-modules and bigraded $(C_p\otimes C_q)$-modules. The compositions $\R \circ {\cal Q}$ and ${\cal Q}\circ \R$ are identities. 
		\end{proof}
		
		\begin{Prop}\label{bimodule1}
			Let $\varphi:\bR^{p-1}\times \bR^{q-1}\to C_{pq}^{00}$ be defined by 
			$$
			\varphi (u_i)= u_i u_p, \ \ \ \varphi(w_i)=w_iw_{q}
			$$
			then $\varphi$ extends to an isomorphism $C_{p-1}\otimes C_{q-1}\simeq (C_p\otimes C_q)^{00}$. 
		\end{Prop}
		\begin{proof}
			We have
			$$
			\varphi(u_i)^2=(u_i u_p)^2=-(u_i)^2 (u_p)^2=-1
			$$
			and
			$$
			\varphi(u_i)\varphi(u_j)=(u_i u_p)( u_j u_p)= -u_i u_j (u_p)^2=u_j u_i (u_p)^2=-(u_j u_p)(u_i u_p)=-\varphi(u_j)\varphi(u_i).
			$$
			Analogically,
			$$
			\varphi(w_i)^2=-1 \ \ \ \ \ \ \text{and} \ \ \ \ \ \ \varphi(w_i)\varphi(w_j) = -\varphi(w_j)\varphi(w_i).
			$$
			Finally,
			$$
			\varphi(u_i)\varphi(w_j)=u_iu_pw_jw_q=w_jw_qu_iu_p=\varphi(w_i)\varphi(u_i).
			$$
			Thus, the element $\varphi(u_1),\ldots \varphi(u_{p}),\varphi_(w_1),\ldots, \varphi(w_q)$ satisfies the relations on generators of the algebra $C_p\otimes C_q$. Therefore, $\varphi$ extends to an isomorphism $C_{p-1}\otimes C_{q-1}\simeq (C_p\otimes C_q)^{00}$. 
		\end{proof}
	Since
 $$
 \bC\otimes \bC \simeq\bC^{ 2}, \ \ \ \bC\otimes \bH \simeq\bC(2), \ \ \ \bH\otimes \bH\simeq\bR(4),
 $$
 for any $n_1,n_2,m_1,m_2 \in \N$ we have 
 \begin{equation}\label{tensorproduct}
 \begin{split}
 \bC(n_1)^{ m_1}\otimes \bC(n_2)^{ m_2} \simeq\bC(n_1n_2)^{2m_1m_2}, \\
 \bC(n_1)^{ m_1}\otimes \bH(n_2)^{ m_1} \simeq\bC(2n_1n_2)^{ m_1m_2}, \\ \bH(n_1)^{ m_1}\otimes \bH(n_2)^{ m_2}\simeq\bR(4n_1n_2)^{ m_1m_2}.
 \end{split}
 \end{equation}
 Hence, we have the multiplication table for $C_p\otimes C_q$.
 
 \medskip
 \begin{tabular}{|c|c|c|c|c|c|c|c|c|}
\hline \ & $C_0$ &$C_1$ & $C_2$ & $C_3$ & $C_4$ & $C_5$ & $C_6$ & $C_7$ \\
\hline $C_0$ & $\bR$& $\bC$ & $\bH$ & $\bH^{ 2}$ & $\bH(2)$ & $\bC(4)$ & $\bR(8)$ & $\bR(8)^{ 2}$ \\
\hline
$C_1$ & $\bC$ & $\bC^{ 2}$ & $\bC(2)$ & $\bC(2)^{ 2}$ & $\bC(4)$ & $\bC(4)^{ 2}$ & $\bC(8)$ & $\bC(8)^{ 2}$ \\
\hline
$C_2$ & $\bH$ & $\bC(2)$ & $\bR(4)$ & $\bR(4)^{ 2}$ & $\bR(8)$ & $\bC(8)$ & $\bH(8)$ & $\bH(8)^{ 2}$ \\
\hline
$C_3$ & $\bH^{ 2}$ & $\bC(2)^{ 2}$ & $\bR(4)^{ 2}$ & $\bR(4)^{ 4}$ & $\bR(8)^{ 2}$ & $\bC(8)^{ 2}$ & $\bH(8)^{ 2}$ & $\bH(8)^{ 4}$ \\
\hline
$C_4$ & $\bH(2)$ & $\bC(4)$ & $\bR(8)$ & $\bR(8)^{ 2}$ & $\bR(16)$ & $\bC(16)$ & $\bH(16)$ & $\bH(16)^{ 2}$ \\
\hline
$C_5$ & $\bC(4)$ & $\bC(4)^{ 2}$ & $\bC(8)$ & $\bC(8)^{ 2}$ & $\bC(16)$ & $\bC(16)^{ 2}$ & $\bC(32)$ & $\bC(32)^{ 2}$ \\
\hline
$C_6$ & $\bR(8)$ & $\bC(8)$ & $\bH(8)$ & $\bH(8)^{ 2}$ & $\bH(16)$ & $\bC(32)$ & $\bR(64)$ & $\bR(64)^{ 2}$ \\
\hline
$C_7$ & $\bR(8)^{ 2}$ & $\bC(8)^{ 2}$ & $\bH(8)^{ 2}$ & $\bH(8)^{ 4}$ & $\bH(16)^{ 2}$ & $\bC(32)^{ 2}$ & $\bR(64)^{ 2}$ & $\bR(64)^{ 4}$ \\
\hline
 \end{tabular}
 \medskip
 \\
 Moreover, since $C_{p+8}\simeq\bR(16)\otimes C_p$, we have 
 $$
 C_{p+8}\otimes C_{q}\simeq C_p\otimes C_{q+8}\simeq\bR(16)\otimes C_p\otimes C_q.
 $$

 \subsection{Classification of \texorpdfstring{$\bZ/2\bZ$}{}-bigraded \texorpdfstring{$C_p\otimes C_q$}{}-modules}

 Let $\bA\in \{\bR,\bC,\bH\}$. Then the real dimension of an indecomposable $\bA(n)^{ m}$-module equals $m\dim_\bR \bA$. Hence, we have the following table for the real dimension $a_{p,q}$ of an indecomposable ungraded $C_p\otimes C_q$-module.

\medskip
 \begin{tabular}{|c|c|c|c|c|c|c|c|c|}
\hline \ & $C_0$ &$C_1$ & $C_2$ & $C_3$ & $C_4$ & $C_5$ & $C_6$ & $C_7$ \\
\hline $C_0$ & $1$& $2$ & $4$ & $4$ & $8$ & 8 & 8 & 8 \\
\hline
$C_1$ & 2 & 2 & 4 & 4 & 8 & 8 & 16 & 16 \\
\hline
$C_2$ & 4 & 4 & 4 & 4 & 8 & 16 & 32 & 32 \\
\hline
$C_3$ & 4 & 4 & 4 & 4 & 8 & 16 & 32 & 32 \\
\hline
$C_4$ & 8 & 8 & 8 & 8 & 16 & 32 & 64 & 64 \\
\hline
$C_5$ & 8 & 8 & 16 & 16 & 32 & 32 & 64 & 64 \\
\hline
$C_6$ & 8 & 16 & 32 & 32 & 64 & 64 & 64 & 64 \\
\hline
$C_7$ & 8 & 16 & 32 & 32 & 64 & 64 & 64 & 64 \\
\hline
\end{tabular} 
\medskip
\\
 Moreover, 
 $$
 a_{p+8,q}=a_{p,q+8}=16a_{p,q}.
 $$

 According to Propositions \ref{bimodule} and \ref{bimodule1}, we have an analogical table for $b_{p,q}=\dim_\bR S_{p,q}^{00}$, where ${S_{p,q}^{00}\oplus S_{p,q}^{01}\oplus S_{p,q}^{10}\oplus S_{p,q}^{11}}$ is an indecomposable bigraded $C_p\otimes C_q$-module. 

 \medskip
 \begin{tabular}{|c|c|c|c|c|c|c|c|c|}
\hline \ & $C_1$ &$C_2$ & $C_3$ & $C_4$ & $C_5$ & $C_6$ & $C_7$ & $C_8$ \\
\hline $C_1$ & $1$& $2$ & $4$ & $4$ & $8$ & 8 & 8 & 8 \\
\hline
$C_2$ & 2 & 2 & 4 & 4 & 8 & 8 & 16 & 16 \\
\hline
$C_3$ & 4 & 4 & 4 & 4 & 8 & 16 & 32 & 32 \\
\hline
$C_4$ & 4 & 4 & 4 & 4 & 8 & 16 & 32 & 32 \\
\hline
$C_5$ & 8 & 8 & 8 & 8 & 16 & 32 & 64 & 64 \\
\hline
$C_6$ & 8 & 8 & 16 & 16 & 32 & 32 & 64 & 64 \\
\hline
$C_7$ & 8 & 16 & 32 & 32 & 64 & 64 & 64 & 64 \\
\hline
$C_8$ & 8 & 16 & 32 & 32 & 64 & 64 & 64 & 64 \\
\hline
\end{tabular} 
\medskip
\\
 Moreover, 
 $$
 b_{p+8,q}=b_{p,q+8}=16b_{p,q}.
 $$

 An algebra $\bA(n)^{ m}$ has $m$ non-isomorphic indecomposable modules. Combining it with the multiplication table and Propositions \ref{bimodule}, \ref{bimodule1} we get the following.
\begin{Prop}
 The number of non-isomorphic ungraded indecomposable $C_p\otimes C_q$-modules equals 
 $$
 \begin{cases}
 4 \ \ \ \text{if} \ \ \ p=4k+3, q=4l+3; \\
 2 \ \ \ \text{if} \ \ \ p=4k+1, q=4l+1 \ \ \ \text{or} \ \ \ p=4k+1, q=4l+3\ \text{or} \ p=4k+3, q=4l+1 ;\\
 1 \ \ \ \text{otherwise}.
 \end{cases}
 $$
\end{Prop}
\begin{Prop}
 The number of non-isomorphic bi-graded indecomposable $C_p\otimes C_q$-modules equals 
 $$
 \begin{cases}
 4 \ \ \ \text{if} \ \ \ p=4k, q=4l; \\
 2 \ \ \ \text{if} \ \ \ p=4k+2, q=4l+2 \ \ \ \text{or} \ \ \ p=4k+2, q=4l\ \text{or} \ p=4k+3, q=4l; \\
 1 \ \ \ \text{otherwise}.
 \end{cases}.
 $$
\end{Prop}

 \section{Classification of Clifford (associative quasi)Nil-algebras of rank 4 and nilpotency index 3} \label{sectionfinal} 

 Now, we finish the classification of rank 4 Clifford Nil-algebras (Theorem \ref{index3}). Actually, we give a more general result. We provide the classification of rank 4 associative Clifford quasNil-algebras.

\begin{theorem}\label{associativeQNAclassification} \ 
\begin{itemize}
 \item [a)]

 Any associative Clifford quasiNil-algebra of rank 4 and nilpotency index 3 admits, up to the ani-transposition, one of the following forms: $\cN_{(0)}^4:=\cN_0^4 ({V^{10}, V^{01}, S^{00},S^{10}, S^{01}, S^{11}})$ where ${S^{00}\oplus S^{10}\oplus S^{01}\oplus S^{11}}$ is a bigraded $\Cl(V^{10}\otimes V^{01})$-module; 
		$$	
		\cN_{(1)}^4:=\begin{pmatrix}
			0 & \bR & V & S^{1} \\
			0 & 0 & V & S^1 \\
			0 & 0 & 0 & S^0 \\
			0 & 0 & 0 & 0
		\end{pmatrix}, \ \ \ 
		\cN_{(2)}^4:=\begin{pmatrix}
			0 & V & V & S^{0} \\
			0 & 0 & \bR & S^1 \\
			0 & 0 & 0 & S^1 \\
			0 & 0 & 0 & 0
		\end{pmatrix}
		$$
		where $S^0\oplus S^1$ is a $\Cl(V)$-module and bilinear maps $\bR\times V \to V, \ \bR\times S^1 \to S^1 $ coincide with standard multiplications;
		$$
		\cN_{(3)}^4:=\begin{pmatrix}
			0 & \C & \C & \C^n \\
			0 & 0 & \C & \C^n \\
			0 & 0 & 0 & \C^n \\
			0 & 0 & 0 & 0
		\end{pmatrix}, \ \ 
		\cN_{(4)}^4=\begin{pmatrix}
			0 & W & \bH & \bH^n \\
			0 & 0 & \bH & \bH^n \\
			0 & 0 & 0 & \bH^n \\
			0 & 0 & 0 & 0
		\end{pmatrix}, \ \ 
		\cN_{(5)}^4:=\begin{pmatrix}
			0 & \bH & \bH & \bH^n \\
			0 & 0 & W & \bH^n \\
			0 & 0 & 0 & \bH^n \\
			0 & 0 & 0 & 0
		\end{pmatrix},
		$$
		where $W\subset \bH$ is a vector subspace.
\item[b)] 
 Let $\cN$ be an associative Clifford quasiNil-algebra of rank 4 and nilpotency index 3. Then $\cN$ is not a Clifford Nil-algebra if and only if $\cN$ is equivalent, up to the ani-transposition, to $\cN_{(2)}^4(V,S^0,S^1)$, where $\dim V>1$, or $\cn_{(4)}^4(W,n)$, where $\dim W<4$.
 \end{itemize} 
\end{theorem}

\begin{proof}
a) Combining Theorem \ref{cN-4}, Proposition \ref{associtiveSQNA}, Theorem \ref{DCEclassification}, and Proposition \ref{bimodulealgebra} we get the Classification of Clifford associative quasiNil-algebras.

 A quasiNil-algebra is a Nil-algebra if and only if the dual quasiNil-algebra is a Nil-algebra. The associative quasiNil-algebra $\cN$ is a Nil-algebra if and only if
 			it satisfies the Vinberg conditions: for any $a_{24} {\in} \cn_{24}$, $a_{34} {\in} \cn_{34}$ 
	$$
 \text{if} \ \ \ \ \langle a_{24}, \cn_{23} \cdot a_{3 4}\rangle {=}\{0\} \ \ \ \
 \text{then} \ \ \ \
 		\langle \cn_{12} \cdot a_{24}, \cn_{13} \cdot a_{34} \rangle{=} \{0\}\ 
 $$
 (Remark \ref{rank4vinbercondition}).
 
 First, suppose
 \begin{equation}\label{first}
 \cN = 
 \begin{pmatrix}
 0 & V & V & S^{0} \\
			0 & 0 & \bR & S^1 \\
			0 & 0 & 0 & S^1 \\
 0 & 0 & 0 & 0
 \end{pmatrix}
 \end{equation}
 
 is a Nil-algebra with $\dim V>1$. According the Vinberg condition, s that for any orthogonal $s_1^1,s_2^2 \in S^1$, we have 
 $$
 \langle Vs_1^1, V s_2^1\rangle=0,
 $$
 but it is not true for 
 $$
 s_1^1=v_1 s^0, \ \ \ \ \ s_2^1=v_2 s^0,
 $$
 where $v_1,v_2\in V$ is orthogonal vectors and $s^0\in S^0$. Thus, there are not Nil-algebras of the form \eqref{first} with $\dim V>1$.

 Second, suppose 
 \begin{equation} \label{second}
 \cN=
 \begin{pmatrix}
			0 & \bH & \bH & \bH^n \\
			0 & 0 & W & \bH^n \\
			0 & 0 & 0 & \bH^n \\
			0 & 0 & 0 & 0
		\end{pmatrix},
 \end{equation}
 is Nil-algebra with $\dim W<4$. According to the the Vinberg condition means that if $q_1, q_2\in \bH$ satisfies the condition
 $$
 \langle wq_1, q_2\rangle=0,
 $$
 then
 $$
 \langle \bH q_1, \bH q_2\rangle=0,
 $$
 $s_1^1,s_2^1 \in S^1$, we have 
 $$
 \langle Vs_1^1, V s_2^1\rangle=0.
 $$
 but it is not true for 
 $$
 q_2=uq_2,
 $$
 where $u\in W^\perp\subset \bH$. Thus, there are not Nil-algebras of the form \eqref{second} with $\dim W<4$.

 Last, suppose a quasiNil-algebra $\cN$ admits one of the following forms 
 
 $$
		\begin{pmatrix}
			0 & V^{10} & S^{10} & S^{11} \\
			0 & 0 & S^{00} & S^{01} \\
			0 & 0 & 0 & V^{01} \\
			0 & 0 & 0 & 0
		\end{pmatrix}, \ \
 \begin{pmatrix}
			0 & S^0 & S^1 & S^{1} \\
			0 & 0 & V & V \\
			0 & 0 & 0 & \bR \\
			0 & 0 & 0 & 0
		\end{pmatrix}, \ \ 
 \begin{pmatrix}
			0 & \C^n & \C^n& \C^n \\
			0 & 0 & \C & \C \\
			0 & 0 & 0 & \C \\
			0 & 0 & 0 & 0
		\end{pmatrix}, \ \ 
 \begin{pmatrix}
			0 & \bH^n & \bH^n& \bH^n \\
			0 & 0 & \bH & \bH \\
			0 & 0 & 0 & \bH\\
			0 & 0 & 0 & 0
		\end{pmatrix}.
		$$
 Then for any nonzero $a_{34}\in \cN_{34}$ the multiplication map
 $$
 R_{a_{34}} : \cN_{23} \to \cN_{24}
 $$
 is an isomorphism. Therefore, for
 any nonzero $ a_{24}\in \cN_{24}, a_{34}\in \cN_{34}$, 
 $$
 \langle a_{24}, \cn_{23} \cdot a_{3 4}\rangle {=} \ \langle a_{24}, \cn_{24}\rangle \ne \{0\}.
 $$
 Thus, the Vinberg condition is satisfied and $\cN$ is Nil-algebra.
\end{proof}

{\bf Acknowledgments:} The work was supported by the Theoretical Physics and Mathematics
Advancement Foundation «BASIS».

	\printbibliography

@article{atiyah1964clifford,
  title={Clifford modules},
  author={Atiyah, Michael F and Bott, Raoul and Shapiro, Arnold},
  journal={Topology},
  volume={3},
  pages={3--38},
  year={1964},
  publisher={Pergamon}
}

@article{andersson2004wishart,
  title={Wishart distributions on homogeneous cones},
  author={Andersson, Steen A and Wojnar, G Gerard},
  journal={Journal of Theoretical Probability},
  volume={17},
  pages={781--818},
  year={2004},
  publisher={Springer}
}

@article{koszul1965varietes,
  title={Vari{\'e}t{\'e}s localement plates et convexit{\'e}},
  author={Koszul, Jean-Louis},
  year={1965}
}

@article{vinberg1963theory,
  title={The theory of convex homogeneous cones},
  author={Vinberg, Ernest B},
  journal={Trans. Moscow Math. Soc.},
  volume={12},
  pages={340--403},
  year={1963}
}

@book{shima2007geometry,
  title={The geometry of Hessian structures},
  author={Shima, Hirohiko},
  year={2007},
  publisher={World Scientific}
}

@article{shima1977symmetric,
  title={Symmetric spaces with invariant locally Hessian structures},
  author={Shima, Hirohiko},
  journal={Journal of the Mathematical Society of Japan},
  volume={29},
  number={3},
  pages={581--589},
  year={1977},
  publisher={The Mathematical Society of Japan}
}

@inproceedings{shima1980homogeneous,
  title={Homogeneous Hessian manifolds},
  author={Shima, Hirohiko},
  booktitle={Annales de l'institut Fourier},
  volume={30},
  number={3},
  pages={91--128},
  year={1980}
}

@article{alekseevsky1997classification,
  title={Classification of N-(super)-extended Poincar{\'e} algebras and bilinear invariants of the spinor representation of Spin (p, q)},
  author={Alekseevsky, Dmitry V and Cort{\'e}s, Vicente},
  journal={Communications in mathematical physics},
  volume={183},
  number={3},
  pages={477--510},
  year={1997},
  publisher={Springer}
}

@article{alekseevsky2021special,
  title={Special Vinberg cones},
  author={Alekseevsky, DV and Cort{\'e}s, Vicente},
  journal={Transformation Groups},
  volume={26},
  pages={377--402},
  year={2021},
  publisher={Springer}
}

@article{hurwitz1898ober,
  title={Ober die Komposition der quadratischen Formen yon beliebig vielen Variablen},
  author={Hurwitz, A},
  journal={Mathematische Werke},
  volume={2},
  year={1898}
}

@article{lee1948theoreme,
  title={Sur le th{\'e}or{\`e}me de Hurwitz-Radon pour la composition des formes quadratiques},
  author={Lee, HC},
  journal={Commentarii Mathematici Helvetici},
  volume={21},
  number={1},
  pages={261--269},
  year={1948},
  publisher={Springer}
}

@article{dubisch1946composition,
  title={Composition of quadratic forms},
  author={Dubisch, Roy},
  journal={Annals of Mathematics},
  volume={47},
  number={3},
  pages={510--527},
  year={1946},
  publisher={JSTOR}
}

@book{lawson2016spin,
  title={Spin Geometry (PMS-38), Volume 38},
  author={Lawson, H Blaine and Michelsohn, Marie-Louise},
  volume={20},
  year={2016},
  publisher={Princeton university press}
}

@article{de1992special,
  title={Special geometry, cubic polynomials and homogeneous quaternionic spaces},
  author={de Wit, Bernard and Van Proeyen, Antoine},
  journal={Communications in Mathematical Physics},
  volume={149},
  pages={307--333},
  year={1992},
  publisher={Springer}
}

@article{cecotti1989geometry,
  title={Geometry of type II superstrings and the moduli of superconformal field theories},
  author={Cecotti, Sergio and Ferrara, Sergio and Girardello, Luciano},
  journal={International Journal of Modern Physics A},
  volume={4},
  number={10},
  pages={2475--2529},
  year={1989},
  publisher={World Scientific}
}

@article{hurwitz1922komposition,
  title={{\"U}ber die Komposition der quadratischen Formen},
  author={Hurwitz, Adolf},
  journal={Mathematische Annalen},
  volume={88},
  number={1},
  pages={1--25},
  year={1922},
  publisher={Springer}
}

@article{sasaki1980hyperbolic,
  title={Hyperbolic affine hyperspheres},
  author={Sasaki, Takeshi},
  journal={Nagoya Mathematical Journal},
  volume={77},
  pages={107--123},
  year={1980},
  publisher={Cambridge University Press}
}

@article{gigena1981conjecture,
  title={On a conjecture by E. Calabi},
  author={Gigena, Salvador},
  journal={Geometriae Dedicata},
  volume={11},
  pages={387--396},
  year={1981},
  publisher={Springer}
}

@inproceedings{calabi1972complete,
  title={Complete affine hyperspheres I},
  author={Calabi, Eugenio},
  booktitle={Symposia Mathematica},
  volume={10},
  pages={19--38},
  year={1972},
  organization={Academic Press London}
}

@article{cheng1986complete,
  title={Complete affine hypersurfaces. Part I. The completeness of affine metrics},
  author={Cheng, Shiu-Yuen and Yau, Shing-Tung},
  journal={Communications on Pure and Applied Mathematics},
  volume={39},
  number={6},
  pages={839--866},
  year={1986},
  publisher={Wiley Online Library}
}

@article{blaschke1923vorlesungen,
  title={Vorlesungen {\"U}ber Differentialgeometrie und Geometrische Grundlagen von Einsteins Relativit{\"a}tstheorie: II Affine Differentialgeometrie},
  author={Blaschke, Wilhelm and Reidemeister, Kurt},
  year={1923},
  publisher={Springer}
}

@article{loftin2008survey,
  title={Survey on affine spheres},
  author={Loftin, John},
  journal={arXiv preprint arXiv:0809.1186},
  year={2008}
}

@article{burde2006left,
  title={Left-symmetric algebras, or pre-Lie algebras in geometry and physics},
  author={Burde, Dietrich},
  journal={Central European Journal of Mathematics},
  volume={4},
  pages={323--357},
  year={2006},
  publisher={Springer}
}

@article{manchon2011short,
  title={A short survey on pre-Lie algebras},
  author={Manchon, Dominique},
  journal={Noncommutative geometry and physics: renormalisation, motives, index theory},
  pages={89--102},
  year={2011},
  publisher={Eur. Math. Soc. Z{\"u}rich}
}

@article{bai2021introduction,
  title={An introduction to pre-Lie algebras},
  author={Bai, Chengming},
  journal={Algebra and Applications},
  volume={1},
  pages={245--273},
  year={2021}
}

@inproceedings{пятецкий1957оценке,
  title={Об оценке размерности пространства автоморфных форм для некоторых типов дискретных групп},
  author={Пятецкий-Шапиро, Илья},
  booktitle={Доклады Академии наук},
  volume={113},
  number={5},
  pages={980--983},
  year={1957},
  organization={Российская академия наук}
}

@book{пятецкий1961геометрия,
  title={Геометрия классических областей и теория автоморфных функций},
  author={Пятецкий-Шапиро, Илья},
  year={1961},
  publisher={Гос. изд-во физико математической лит-ры}
}

@article{riesz1949integrale,
  title={L'int{\'e}grale de Riemann-Liouville et le probl{\`e}me de Cauchy},
  journal={Acta Math},
  volume={81},
  author={Riesz, Marcel},
  year={1949},
 pages={1--221}
}

@article{gindikin1964analysis,
  title={Analysis in homogeneous domains},
  author={Gindikin, Simon G},
  journal={Russian Mathematical Surveys},
  volume={19},
  number={4},
  pages={1},
  year={1964},
  publisher={IOP Publishing}
}

@article{ishi2000positive,
  title={Positive Riesz distributions on homogeneous cones},
  author={Ishi, Hideyuki},
  journal={Journal of the Mathematical Society of Japan},
  volume={52},
  number={1},
  pages={161--186},
  year={2000},
  publisher={The Mathematical Society of Japan}
}

@article{hassairi2005beta,
  title={Beta-Riesz distributions on symmetric cones},
  author={Hassairi, A and Lajmi, S and Zine, R},
  journal={Journal of statistical planning and inference},
  volume={133},
  number={2},
  pages={387--404},
  year={2005},
  publisher={Elsevier}
}

@article{said2017riemannian,
  title={Riemannian Gaussian distributions on the space of symmetric positive definite matrices},
  author={Said, Salem and Bombrun, Lionel and Berthoumieu, Yannick and Manton, Jonathan H},
  journal={IEEE Transactions on Information Theory},
  volume={63},
  number={4},
  pages={2153--2170},
  year={2017},
  publisher={IEEE}
}

@incollection{said2022gaussian,
  title={Gaussian distributions on Riemannian symmetric spaces of nonpositive curvature},
  author={Said, Salem and Mostajeran, Cyrus and Heuveline, Simon},
  booktitle={Handbook of Statistics},
  volume={46},
  pages={357--400},
  year={2022},
  publisher={Elsevier}
}

@article{alekseevsky2009geometric,
  title={Geometric construction of the r-map: from affine special real to special K{\"a}hler manifolds},
  author={Alekseevsky, Dmitri V and Cort{\'e}s, Vicente},
  journal={Communications in Mathematical Physics},
  volume={291},
  number={2},
  pages={579--590},
  year={2009},
  publisher={Springer}
}

@article{alekseevsky2013conification,
  title={Conification of K{\"a}hler and hyper-K{\"a}hler manifolds},
  author={Alekseevsky, Dmitri V and Cort{\'e}s, Vicente and Mohaupt, Thomas},
  journal={Communications in Mathematical Physics},
  volume={324},
  pages={637--655},
  year={2013},
  publisher={Springer}
}

@article{alekseevskii1975structure,
  title={Structure of homogeneous Riemann spaces with zero Ricci curvature},
  author={Alekseevskii, Dmitrii  and Kimel'fel'd, BN},
  journal={Functional Analysis and its Applications},
  volume={9},
  number={2},
  pages={97--102},
  year={1975},
  publisher={Springer}
}

@article{cortes1996homogeneous,
  title={Homogeneous special geometry},
  author={Cort{\'e}s, Vicente},
  journal={Transformation Groups},
  volume={1},
  number={4},
  pages={337--373},
  year={1996},
  publisher={Springer}
}

@article{alekseevsky2021special1,
  title={Special Vinberg cones and the entropy of BPS extremal black holes},
  author={Alekseevsky, Dmitri V and Marrani, Alessio and Spiro, Andrea},
  journal={Journal of High Energy Physics},
  volume={2021},
  number={11},
  pages={1--35},
  year={2021},
  publisher={Springer}
}

@article{alekseevsky2023erratum,
  title={Erratum to: Special Vinberg cones and the entropy of BPS extremal black holes},
  author={Alekseevsky, Dmitri and Marrani, Alessio and Spiro, Andrea},
  journal={Journal of High Energy Physics},
  volume={2023},
  number={2107.06797},
  pages={1--3},
  year={2023},
  publisher={Springer}
}

@article{cortes2017class,
  title={A class of cubic hypersurfaces and quaternionic K$\backslash$" ahler manifolds of co-homogeneity one},
  author={Cort{\'e}s, Vicente and Dyckmanns, Malte and J{\"u}ngling, Michel and Lindemann, David},
  journal={arXiv preprint arXiv:1701.07882},
  year={2017}
}

@article{bohm2022homogeneous,
  title={Homogeneous Einstein metrics on Euclidean spaces are Einstein solvmanifolds},
  author={B{\"o}hm, Christoph and Lafuente, Ramiro A},
  journal={Geometry \& Topology},
  volume={26},
  number={2},
  pages={899--936},
  year={2022},
  publisher={Mathematical Sciences Publishers}
}

@inproceedings{nemirovski2006advances,
  title={Advances in convex optimization: conic programming},
  author={Nemirovski, Arkadi},
  booktitle={International Congress of Mathematicians},
  volume={1},
  pages={413--444},
  year={2006}
}

@book{nesterov1994interior,
  title={Interior-point polynomial algorithms in convex programming},
  author={Nesterov, Yurii and Nemirovskii, Arkadii},
  year={1994},
  publisher={SIAM}
}

@book{ben2001lectures,
  title={Lectures on modern convex optimization: analysis, algorithms, and engineering applications},
  author={Ben-Tal, Aharon and Nemirovski, Arkadi},
  year={2001},
  publisher={SIAM}
}

@article{shmakova1997calabi,
  title={Calabi-Yau black holes},
  author={Shmakova, Marina},
  journal={Physical Review D},
  volume={56},
  number={2},
  pages={R540},
  year={1997},
  publisher={APS}
}

@article{matrix,
  title={Homogeneous cones associated with Clifford modules},
  author={Alekseevsky, Dmitri V and Cort{\'e}s, Vicente},
  journal={Matrix},
  volume={Special   issue  dedicated     to  the memory of Joseph  Wolf},
  number={to appear},
  year={2025},
  publisher={APS}
}

@book{shapiro2011compositions,
  title={Compositions of quadratic forms},
  author={Shapiro, Daniel B},
  volume={33},
  year={2011},
  publisher={Walter de Gruyter}
}

@book{schneider2022convex,
  title={Convex Cones: Geometry and Probability},
  author={Schneider, Rolf},
  volume={2319},
  year={2022},
  publisher={Springer}
}

@article{wishart1928generalised,
  title={The generalised product moment distribution in samples from a normal multivariate population},
  author={Wishart, John},
  journal={Biometrika},
  volume={20},
  number={1/2},
  pages={32--52},
  year={1928},
  publisher={JSTOR}
}

@article{morris2014natural,
  title={Natural exponential families},
  author={Morris, Carl N},
  journal={Wiley StatsRef: Statistics Reference Online},
  year={2014},
  publisher={Wiley Online Library}
}


\end{document}